\documentclass[11pt]{article}

\usepackage{times,amssymb,amsmath,exscale,array,latexsym}
\usepackage{graphicx}
\usepackage{epsfig}
\usepackage{subfig}
\usepackage{color}
\definecolor{marin}{rgb}   {0.,   0.3,   0.7}
\definecolor{rouge}{rgb}   {0.8,   0.,   0.}
\definecolor{sepia}{rgb}   {0.8,   0.5,   0.}
\usepackage[colorlinks,citecolor=marin,linkcolor=rouge,
            bookmarksopen,
            bookmarksnumbered
           ]{hyperref}

\newcommand{\e}{\ensuremath{\mathrm{e}}}

\newcommand{\psis}{{\mathcal S}}
\renewcommand{\r}{\mathbb{R}}

\newtheorem{lemma}{Lemma}[section]

\newtheorem{remark}[lemma]{Remark}

\numberwithin{equation}{section}

\newcommand{\QED}{\mbox{}\hfill \raisebox{-0.2pt}{\rule{5.6pt}{6pt}\rule{0pt}{0pt}}
          \medskip\par}

\begin{document}
\title{Runge--Kutta--Nystr\"om symplectic splitting methods of order 8
}

\author{S. Blanes$^{1}$, F. Casas$^{2}$, A. Escorihuela-Tom\`as$^{3}$ \\[2ex]
$^{1}$ {\small\it Universitat Polit\`ecnica de Val\`encia, Instituto de Matem\'atica Multidisciplinar, 46022-Valencia, Spain}\\{
\small\it email: serblaza@imm.upv.es}\\[1ex]
$^{2}$ {\small\it Departament de Matem\`atiques and IMAC, Universitat Jaume I, 12071-Castell\'on, Spain}\\{
\small\it email: Fernando.Casas@mat.uji.es}\\[1ex]
$^{3}$ {\small\it Departament de Matem\`atiques, Universitat Jaume I, 12071-Castell\'on, Spain}\\{
\small\it email: alescori@uji.es}\\[1ex]
}


\date{May 3, 2022}

\maketitle

\begin{abstract}

Different families of  Runge--Kutta--Nystr\"om (RKN) symplectic splitting methods of order 8 are presented for second-order systems of ordinary differential equations and are tested on 
 numerical examples. They show a better efficiency than state-of-the-art symmetric compositions of 2nd-order symmetric schemes and
RKN splitting methods of orders 4 and 6 for medium to high accuracy. For some particular examples, they are even more efficient than
extrapolation methods for high accuracies and integrations over relatively short time intervals.



\end{abstract}\bigskip


\section{Introduction}
\label{sec.1}

Second-order systems of ordinary differential equations (ODEs) of the form
\begin{equation}  \label{rkn.1}
  \ddot{y} \equiv \frac{d^2 y}{dt^2} =  g(y),
\end{equation}  
where $y \in \mathbb{R}^d$ and $g: \mathbb{R}^d \longrightarrow  \mathbb{R}^d$, appear very often in applications, so that special numerical
integrators have been designed for them, such as the Runge--Kutta--Nystr\"om (RKN) class of methods. As is well known, if one introduces the
new variables $x = (y,v = \dot{y})$ and the maps
\begin{equation} \label{rkn.2}
  f_a(x) = f_a(y,v) = (v,0), \qquad f_b(x) = f_b(y,v) = (0, g(y)),
\end{equation}
then eq. (\ref{rkn.1}) is equivalent to
\begin{equation} \label{rkn.3}
  \dot{x} = f_a(x) + f_b(x)
\end{equation}
and moreover each subsystem $\dot{x} = f_i(x)$, $i=a,b$, is explicitly integrable, with exact flow
\[
 \varphi_t^{[a]}(y,v) = (y + t v, v) \qquad \mbox{ and } \qquad  \varphi_t^{[b]}(y,v) = (y, v+ t g(y)),
\]
respectively.   An important class of problems leading to equations of the form (\ref{rkn.1}) corresponds to Hamiltonian dynamical systems of the form
\begin{equation} \label{ham1}
  H(q,p) = \frac{1}{2} p^T M^{-1} p + V(q),
\end{equation}
where $q$ and $p$ denote coordinates and momenta, respectively, $M$ is a symmetric positive definite 
square constant matrix and $V(q)$ is the potential. Then,
the corresponding equations of motion can be written as (\ref{rkn.1}) with $y = q$, $v = \dot{y} = M^{-1} p$ and $g(y) = - {M}^{-1}\nabla V(q)$.

Splitting methods constitute a natural option for integrating numerically the initial value problem defined by (\ref{rkn.3}).
These are schemes of the form
\begin{equation} \label{splitting1}
  \psi_h = \varphi_{h a_s}^{[a]} \circ \varphi_{h b_s}^{[b]} \circ \cdots \circ  \varphi_{h a_1}^{[a]} \circ \varphi_{h b_1}^{[b]}, 
\end{equation}
where the coefficients $a_j$, $b_j$ are conveniently chosen so as to achieve high order approximations to the exact flow of (\ref{rkn.3}), namely
$\varphi_h(x) = \psi_h(x) + \mathcal{O}(h^{r+1})$ for a given order $r$ and step size $h$. Familiar examples of splitting methods
are the so-called Strang/leapfrog/St\"ormer--Verlet
second order schemes:
\begin{equation} \label{strang.a}
  \mathcal{S}_h^{[2]} = \varphi_{h/2}^{[a]} \circ  \varphi_{h}^{[b]} \circ  \varphi_{h/2}^{[a]}, 
\end{equation}
and
\begin{equation} \label{strang.b}  
  \mathcal{S}_h^{[2]} = \varphi_{h/2}^{[b]} \circ  \varphi_{h}^{[a]} \circ  \varphi_{h/2}^{[b]}.
\end{equation}
In fact, efficient schemes 
of this class up to order $r=6$ 
have been designed along the years (see e.g. \cite{blanes08sac} and references therein). 
In addition, they preserve qualitative properties of the continuous system  and
show a very good behavior with respect to the propagation of errors, especially for long time integrations \cite{hairer06gni}.

There are situations, however, when even higher-order numerical approximations ($r=8,10,\ldots$) are
required, for instance in problems arising in astrodynamics. In that case, although generic splitting methods exist, they involve such a large number
of elementary flows $\varphi_h^{[a]}$,  $\varphi_h^{[b]}$, that are not competitive with other integrators. This is so due to the exponential growth with the order $r$ 
of the required number of
conditions to be satisfied to achieve that order \cite{mclachlan02sm}. For this reason, palindromic compositions of the form
\begin{equation} \label{eq.2.1.2}
   \psis_{\alpha_{m} h}^{[2]}\circ \,  \psis_{\alpha_{m-1} h}^{[2]} \circ
\cdots\circ
 \psis_{\alpha_{2} h}^{[2]}  \circ \, \psis_{\alpha_{1}h}^{[2]}  \quad\qquad \mbox{ with } \quad 
(\alpha_1,\ldots,\alpha_{m}) \in \r^{m}
\end{equation}
and $\alpha_{m+1-i} = \alpha_i$, have been considered instead for order $r > 6$. In practice, schemes (\ref{eq.2.1.2}) are the most realistic option
when one is interested in integrating (\ref{rkn.3}) with high-order ($r=8,10,\ldots$) splitting methods.

It turns out, however, that the special structure of (\ref{rkn.2})-(\ref{rkn.3}) corresponding to the system (\ref{rkn.1}) leads to a reduction in the number of
order conditions when $r > 4$ with respect to the generic problem. This allows one to construct highly efficient 4th- and 6th-order splitting methods especially tailored for this class of
problems which show a better performance than schemes of the family (\ref{eq.2.1.2}) \cite{blanes02psp,omelyan02otc}. 
They can be naturally called 
\emph{RKN splitting methods},
and the question of the existence of eighth-order schemes, more efficient than methods of type (\ref{eq.2.1.2}), formulated some 25 years ago
\cite[p. 153]{mclachlan96aso}, still
remains unanswered, no doubt due to the technical difficulties involved. 

It is our purpose in this note to
 present new RKN splitting methods of order 8 that provide higher efficiency than state-of-the-art composition methods
(\ref{eq.2.1.2}) on a variety of examples arising in physical applications. They should then 
be considered as the natural option when one is interested in integrating numerically
problems of the form (\ref{rkn.2})-(\ref{rkn.3}) with medium to high precision whereas preserving by construction the main qualitative features of the continuous system.

\begin{remark}
It turns out that this class of schemes can also be used to solve the slightly more general problem
\begin{equation}  \label{rkn.1b}
  \ddot{y} = \alpha \dot{y} + \beta y + g(t,y),
\end{equation}  
where $\alpha, \beta\in \mathbb{R}^{d\times d}$ are constant: by taking time $t$ as a new coordinate and considering $x = (y, v, t)$, it is clear that  equation
(\ref{rkn.1b}) can be again expressed as (\ref{rkn.3}), this time with
\begin{equation} \label{rkn.2b}
  f_a(x) = f_a(y,v,t) = (v,\alpha v + \beta y,1), \qquad f_b(x) = f_b(y,v,t) = (0, g(t,y),0),
\end{equation}
and each sub-system being explicitly integrable.
\end{remark}

\section{Order conditions}

As shown e.g. in \cite{blanes16aci}, to each integrator (\ref{splitting1}) one can associate  a series $\Psi(h)$ of differential operators given by
\begin{equation} \label{ser1}
  \Psi(h) = \exp(h b_1 F_b) \, \exp(h a_1 F_a) \cdots \exp(h b_s F_b) \, \exp(h a_s F_a),
\end{equation}
where $F_a$, $F_b$ are the Lie derivatives corresponding to $f_a$ and $f_b$, respectively \cite{arnold89mmo}: 
for each smooth function $g: \mathbb{R}^d \longrightarrow
\mathbb{R}^d$ and $x \in \mathbb{R}^d$ one has
\begin{equation} \label{der.Lie}
  F_a \, g(x) = f_a(x) \cdot \nabla g(x), \qquad\qquad  F_b \, g(x) = f_b(x) \cdot \nabla g(x),
\end{equation}
so that, for the whole integrator, $g(\psi_h(x)) = \Psi(h) g(x)$. For $g(x)=(g_1(x),\ldots,g_d(x))$, we denote
\[
  f(x) \cdot\nabla g(x) \equiv (f(x) \cdot\nabla g_1(x),\ldots,f(x) \cdot\nabla g_d(x))
\]
in eq. (\ref{der.Lie}).  
The main advantage of using the series $\Psi(h)$ for representing the method $\psi_h$ is
that one can formally apply the Baker--Campbell--Hausdorff formula \cite{varadarajan84lgl} and express $\Psi(h)$ as only one exponential,
\begin{equation} \label{ser2}
  \Psi(h) = \exp(F(h)), \qquad \mbox{ with } \qquad F(h) = \sum_{j \ge 1} h^j F_j,
\end{equation}
and each $F_j$ is a linear combination of nested commutators involving $j$ operators $F_a$ and $F_b$ whose coefficients are polynomials of degree $j$
in the coefficients $a_i$, $b_i$. A method of order $r$ requires that $F_1 = F_a + F_b$ for consistency, and $F_j = 0$ for $1 < j \le r$. These constraints
in turn lead to a set of polynomial equations to be satisfied by the coefficients of the splitting method. The number $n_r$ of such order conditions at each $r$
is collected in Table \ref{table.1} \cite{mclachlan02sm}. For comparison, we also include the number $s_r$ of order conditions for compositions of the form
(\ref{eq.2.1.2})

\begin{table}[!h]
 \begin{center}
  \begin{tabular}{|c|cccccccccc|} \hline
  $r$    & 1 & 2 & 3 & 4 & 5 & 6 & 7  &  8 & 9 & 10  \\ \hline
{$s_r$}    & 1 & 0 & 1 & 1 & 2 & 2 & 4 & 5 & 8 & 11  \\ 
$n_r$    & 2 & 1 & 2 & 3 & 6 & 9 & 18 & 30 & 56 & 99  \\ 
$\ell_r$ & 2 & 1 & 2 & 2 & 4 & 5 & 10 & 14 & 25 & 39  \\ \hline
\end{tabular}
\end{center}
 \caption{{\small Number of independent order conditions (at order $r$) of compositions of symmetric second order methods of the form
 (\ref{eq.2.1.2}), $s_r$, of splitting methods in the general case, $n_r$, and in the RKN case, $\ell_r$.}}
  \label{table.1}
\end{table}

As is well known, if the composition (\ref{splitting1}) is left-right palindromic, then all the order conditions at even order
are automatically satisfied and the method is time-symmetric. 
For systems of the form (\ref{rkn.2})-(\ref{rkn.3}), the flow $\varphi_h^{[b]}$ is typically the most expensive part to evaluate (for the Hamiltonian (\ref{ham1}), it corresponds
 essentially to the force $\nabla V(q)$). It makes sense, then, to characterize a given splitting method according to the number of flows
 $\varphi_h^{[b]}$ involved. This is called the \emph{number of stages} of the method. Notice that, if the Strang splitting is used as
 the scheme $\mathcal{S}_h^{[2]}$ in the composition (\ref{eq.2.1.2}), the number of stages is also $m$.
 
 From Table~\ref{table.1} it is then straightforward to estimate the minimum number of stages to achieve an even order $r=2k$. For the composition \eqref{eq.2.1.2} and the general splitting (\ref{splitting1}) these values are, respectively,
 \[
 S_r=2\sum_{i=1}^ks_{2k-1}-1, \qquad 
 N_r=\sum_{i=1}^kn_{2k-1}-1, \qquad 
\]
and are collected in Table~\ref{table.2} up to $r=2k=10$. Notice that, when counting the number of stages per step, we have used the so--called FSAL (First 
Same As Last) property: the last map in one step can be saved in the following one and does not count for the total number of stages.

\begin{table}[!h]
 \begin{center}
  \begin{tabular}{|c|ccccc|} \hline
  $r$  & 2 & 4 & 6 &  8 & 10  \\ \hline
$S_r$  & 1 & 3 & 7 & 15 & 31   \\ 
$N_r$  & 1 & 3 & 9 & 27 & 83   \\ 
$L_r$  & 1 & 3 & 7 & 17 & 42   \\ \hline
\end{tabular}
\end{center}
 \caption{{\small Minimum number of stages required to achieve order $r=2k$ with symmetric compositions (\ref{eq.2.1.2}), $S_r$, with general splitting
 (\ref{splitting1}), $N_r$, and for RKN splitting methods, $L_r$.}}
  \label{table.2}
\end{table}

The number of order conditions to be solved for each family of methods is, respectively, $(S_r+1)/2$ and $N_r+1$. It is clear that symmetric compositions
(\ref{eq.2.1.2}) require to solve a considerably smaller number of order conditions to achieve high order methods. On the other hand, the space of solutions is significantly
larger in the case of general splitting methods, and consequently also the chance of finding highly efficient schemes within this class. Thus, in particular, 
the general splitting methods of order four and six presented in \cite{blanes02psp} outperform compositions (\ref{eq.2.1.2}) of the same order. At order eight,
however, one has to solve a system of 28 polynomial equations for general splitting methods, and although it seems quite likely that very efficient solutions
exist, to carry out out a thorough analysis constitutes a formidable task.


Notice that for systems of the form (\ref{rkn.2})-(\ref{rkn.3}) one has further restrictions: since
$F_a=v \, \nabla_y$, and $F_b=g(y) \nabla_v$, one has for symmetric methods
\[
  [F_b,[F_a, F_b]] = \tilde g(y) \nabla_v, \qquad \mbox{\rm with} \qquad 
	\tilde g(y)=2 \, \nabla_y g(y)\cdot g(y),
\]
where $[F_a, F_b]=F_aF_b-F_bF_a$, etc. In consequence,
$[F_b,[F_b,[F_a, F_b]]] \equiv 0$, and many terms in
(\ref{ser2}) vanish identically, so that their order conditions can be ignored. This can be seen in the last row of Tables \ref{table.1} and \ref{table.2}, where
 we collect the order conditions $\ell_r$ and the the minimum number of stages, 
 \[
    L_r=\sum_{i=1}^k\ell_{2k-1}-1
 \]  
 up to $r=10$. Notice that, whereas the reduction up to $r=6$ with respect to general splitting methods is only of two equations, 
 for a time-symmetric method of order $r=8$ one has
to solve 18 order conditions (instead of 28). This problem, although more amenable, is still far from trivial. 
In addition, to get significant solutions, the relevant issue here is whether the resulting 8th-order 
RKN splitting schemes are competitive
in terms of the number of flows involved with methods within the class (\ref{eq.2.1.2}). 

\begin{remark}
With respect to the more general system \eqref{rkn.1b}-\eqref{rkn.2b}, one has
\[
  F_a=v \, \nabla_y+(\alpha v + \beta y) \nabla_v+1\cdot \partial_t, \qquad  F_b=g(t,y) \nabla_v,
\]
so that  
\[
  [F_b,[F_a, F_b]] = \tilde g(t,y) \nabla_v, \qquad \mbox{\rm with} \qquad 
	\tilde g(y)=2 \nabla_y g(t,y)\cdot g(t,y)
\]
and therefore $[F_b,[F_b,[F_a, F_b]]] \equiv 0$ also here.
\end{remark}

Before starting a systematic search of solutions to the order conditions, it seems appropriate to make explicit several considerations: 
\begin{enumerate} 
 \item Due to the different qualitative character of the operators $F_a$ and $F_b$, it is clear that the role of $\varphi_h^{[a]}$ and $\varphi_h^{[b]}$ in 
 (\ref{splitting1}) is not interchangeable, and so two different orderings have to be considered. Specifically, we will analyze two types of composition:
\begin{equation} \label{aba}
  \mathcal{A}_s =  \varphi_{h a_{s+1}}^{[a]} \circ \varphi_{h b_s}^{[b]} \circ \varphi_{h a_s}^{[a]} \circ \cdots   \circ \varphi_{h b_1}^{[b]} \circ  \varphi_{h a_1}^{[a]},
\end{equation}
with $a_{s+2-i} = a_i$, $b_{s+1-i} = b_i$, and
\begin{equation} \label{bab}
  \mathcal{B}_s =  \varphi_{h b_{s+1}}^{[b]} \circ \varphi_{h a_{s}}^{[a]} \circ \varphi_{h b_s}^{[b]} \circ  \cdots  \circ \varphi_{h a_1}^{[a]} \circ \varphi_{h b_1}^{[b]},   
\end{equation}
with $b_{s+2-i} = b_i$, $a_{s+1-i} = a_i$. Since for methods (\ref{aba}) and (\ref{bab}) one can always apply the FSAL property,
we say that
 both schemes involve the same number $s$ of stages.
 \item Very often, compositions with a higher number of stages than the minimum required to solve the order conditions are considered in the literature. 
 This is so because, typically, (i) methods with the minimum number of stages show a poor performance, and (ii) the presence of free parameters allows
 one to optimize the schemes according with some appropriate criteria, so that the extra computational cost is compensated by the reduction in the error. 
 Thus, in particular, 8th-order methods within the class (\ref{eq.2.1.2}) with 17, 19 and 21 stages
 exist that are more efficient than schemes with the minimum number $m=15$. Notice in this respect that the minimum number of stages for a RKN splitting method of order 8
  is $s=17$. Although one such method of the form $\mathcal{A}_s$ was proposed in \cite{okunbor94oeo}, the numerical results collected there show no clear improvement
 with respect to the 8th-order method of type (\ref{eq.2.1.2}) with $m=24$ presented in \cite{calvo93hos}.
 \item Given a method $\psi_h$, one may consider a near-to-identity map $\pi_h$ so that the
 integrator  $\hat{\psi}_h = \pi_h^{-1} \circ \psi_h \circ \pi_h$ is more accurate than $\psi_h$, for instance, by increasing its order. In this context,
 $\psi_h$ is called the kernel of the processed method $\hat{\psi}_h$, and $\pi_h$ is the processor or corrector.
 Notice that $N$
 consecutive steps correspond to $\hat{\psi}_h^N = \pi_h^{-1} \circ \psi_h^N \circ \pi_h$, i.e., the cost of applying the processed scheme is basically the
 cost of the kernel. This technique allows one
 to separate the order conditions into two sets: the conditions satisfied by the kernel itself, and those to be verified by the processor. As a result, it is possible
 to construct  high-order RKN splitting methods involving a reduced number of stages in the kernel, although building a particular processor is far
 from trivial. Methods of this class have been presented in \cite{blanes01hor,blanes01nfo}, so that they will not be considered here.
 \item For the initial value problem defined by (\ref{rkn.2})-(\ref{rkn.3}), it is possible to include in the compositions (\ref{aba}) and (\ref{bab}) the flows
 generated by other vector fields lying in the Lie algebra generated by $F_a$ and $F_b$. For instance, one could use the $h$-flow of the vector fields
 $[F_b,[F_a,F_b]]$, $[F_b,[F_b,[F_a,[F_a,F_b]]]]$, and other more general nested commutators \cite{blanes01hor,blanes01nfo}. These give rise to the so-called
 `modified potentials', and allow one to reduce the number of stages (although at the price of an additional computational cost to evaluate the flows). Methods
 of this class with and without processing have been analyzed in particular in \cite{blanes01hor} and \cite{omelyan02otc}. Here, by contrast, we are only
 interested in standard compositions (\ref{aba})-(\ref{bab}).
 
 \end{enumerate}

\section{New methods of order 8}

We next analyze families of schemes (\ref{aba}) and (\ref{bab}) involving $s=17, 18$ and 19 stages, so that one always 
has enough parameters in the compositions to solve the 
order conditions. 
Of course, even with the minimum number of parameters, these order conditions possess a large
 number of real solutions, so that
some criterion has to be adopted to select ``good'' methods. 
As is customary in the literature, and assuming $h$ is sufficiently small and $g$ is sufficiently smooth, we propose
to take the leading term in the asymptotic expansion of the modified vector field associated with the integrator as the main contribution to the truncation error. Without any specific 
assumption on the function $g$, we take this error as $(\sum_{i=1}^{25} k_{9,i}^2)^{1/2}$. Here $k_{9,i}$ are the coefficients of the asymptotic expansion
of the modified vector field at order $h^9$ when it is expressed as a linear combination of the 25 independent nested commutators involving 9 operators
$F_a$ and $F_b$. This corresponds to the subspace of the Lie algebra generated by $F_a$ and $F_b$ with the commutator as the Lie bracket (for more details, see \cite{mclachlan03tae,mclachlan19tla}). 
To take into account the computational cost, we multiply this error by the
number of stages $s$, thus resulting in the following effective error for a method of order 8,
\begin{equation} \label{eferror}
  E_{\mathrm{f}} = s \cdot \left(\sqrt{\sum_{i=1}^{25} k_{9,i}^2}\right)^{1/8},
\end{equation}
which should be minimized by the integrator.
One has to take into account, however, that the expression of $E_{\mathrm{f}}$ depends on the particular basis of nested commutators one is considering
and that we are also assuming that all these commutators contribute in a similar way, something that is not guaranteed to take place in all applications. 
It makes sense, then, to introduce other quantities as possible
estimators of the error committed. In particular, it has been noticed that large coefficients $a_i$, $b_i$ in the splitting method 
usually leads to large truncation errors, since the expressions of  $k_{\ell,j}$  for $\ell \ge 9$ depend on increasingly higher powers of these coefficients. 
For this reason, we also keep track of the quantities
\begin{equation} \label{sizeco}
  \Delta \equiv \sum_{i=1}^s \left( |a_i| + |b_i| \right) \qquad \mbox{ and } \qquad \delta \equiv \max_{i=1}^s \, (|a_i|, |b_i|)
\end{equation}
and eventually discard solutions with large values of $\Delta$ and/or $\delta$.  By following a similar approach as for instance in \cite{blanes02psp,omelyan02otc}, we will select particular schemes with small values of $E_{\mathrm{f}}$, $\Delta$ and $\delta$, and then we will test them on an array of numerical examples to check their
efficiency in practice.

\paragraph{$s=17$ stages.}
In this case one has as many parameters as order conditions, 18 in total. Given the complexity of the problem, it is not possible to solve these nonlinear
equations with a computer algebra system, and so one has to turn to numerical techniques. Specifically, they
are solved with the Python~\cite{python3rm} function \texttt{fsolve} of the
{\it SciPy} library~\cite{scipy}, a wrapper of the classic subroutines {HYBRD} and {HYBRJ} of {MINPACK}~\cite{minpack}. The algorithm is
 based on a modification of
the Powell hybrid method and involves the choice of the correction as a convex combination of the Newton method and scaled gradient directions and the
updating of the Jacobian by the rank-1 method (except at the starting point, where it is approximated by forward differences). 
Since we are not
interested in methods with large values of $\delta$, a uniform
distribution in the interval $[-1,1]$ in each variable was taken to generate about $2 \times 10^6$ initial points to start the procedure, 

When a composition of type $\mathcal{A}_s$ is considered, we have obtained 376 real solutions that cannot be obtained as a
composition of 2nd-order symmetric schemes (\ref{eq.2.1.2}), with parameters 
$E_{\mathrm{f}} \in [2.77, 18.05]$ and $\Delta \in [8.40, 63.05]$, respectively. Among these, we select 
those solutions within the more restricted range
$E_{\mathrm{f}} \in [2.86, 3.45]$ and $\Delta \in [8.42, 19.30]$ and check them on the test problems of sections \ref{num.test1} and \ref{num.test2}. Finally, we have
chosen the scheme whose coefficients are listed in Table \ref{tau:rkna}, and parameters  given in 
Table \ref{tau:ef}. The final values of the coefficients (with 30 digits of accuracy) have been obtained by taking the solution found by \texttt{fsolve} as the
starting point of the function \texttt{FindRoot} of \emph{Mathematica}.
The method can be represented in the compact form
\begin{align}\label{eq:rkna17}
  \mathcal{A}_{17} \, \equiv \, & (a_1, b_1, a_2, b_2, a_3, b_3, a_4, b_4, a_5, b_5, a_6, b_6, a_7, b_7, a_8, b_8,a_9,b_9, \nonumber\\
  \quad \quad &a_9,b_8,a_8, b_7, a_7, b_6, a_6, b_5, a_5, b_4, a_4, b_3, a_3, b_2, a_2, b_1, a_1).
\end{align}

For compositions of type $\mathcal{B}_s$, by applying the same methodology, we have found 149 different solutions out of more than $1.2 \times 10^6$  starting points. We have selected the four solutions in the region $E_{\mathrm{f}} \in [2.80, 3.85]$, $\Delta \in  [7.30, 9.95]$ and finally we take the one whose
coefficients are collected in Table \ref{tau:rknb}. The method thus reads
\begin{align}
  \mathcal{B}_{17} \, \equiv \, & (b_1, a_1, b_2, a_2, b_3, a_3, b_4, a_4, b_5, a_5, b_6, a_6, b_7, a_7, b_8, a_8,b_9,a_9 \nonumber\\
  \quad \quad &b_9,a_8,b_8, a_7, b_7, a_6, b_6, a_5, b_5, a_4, b_4, a_3, b_3, a_2, b_2, a_1, b_1).  
\end{align}

\begin{table}[ht]
  \renewcommand\arraystretch{1.4}
  \begin{center}
    \begin{tabular}{lllllll}
      &&\multicolumn{1}{c}{$E_{\rm{f}}$}&&\multicolumn{1}{c}{$\Delta$} && \multicolumn{1}{c}{$\delta$} \\ 
      \cline{1-7}
      $\mathcal{A}_{17}$&\hspace{1.5cm} &$3.45 $&\hspace{0.4cm} &$8.42$ & \hspace{0.4cm} &$0.5459$ ($|a_9|$) \\
      $\mathcal{A}_{18}$& &$3.65 $ &&$7.42 $ && $0.6406$ ($|a_9|$) \\
      $\mathcal{A}_{19}$& &$2.76 $ &&$5.98 $ && $0.4237$ ($|a_4|$) \\
      $\mathcal{B}_{17}$& &$2.80 $ &&$8.93 $ & & $0.6355$ ($|a_5|$) \\
      $\mathcal{B}_{18}$& &$3.44 $ &&$9.68 $ && $0.9303$ ($|a_4|$) \\
      $\mathcal{B}_{19}$& &$3.41 $ &&$6.94 $ && $0.5238$ ($|a_6|$) \\
      \cline{1-7}
    \end{tabular}
    \caption{{\small Effective error  $E_{\mathrm{f}}$, 1- and $\infty$-norm of the vector of coefficients for different  8th-order RKN splitting methods
    of type $\mathcal{A}_s$ and $\mathcal{B}_s$.}}
    \label{tau:ef}
  \end{center}
\end{table}

\begin{table}[!h]
{\small
  \renewcommand\arraystretch{1.1}
  \begin{center}
    \begin{tabular}{lll}
      &\multicolumn{1}{c}{$a_i$}&\multicolumn{1}{c}{$b_i$} \\ 
      \cline{2-3}
      $\mathcal{A}_{17}$&$a_1=0.0520924343840339006426037968353$ &$b_1=0.145850304812644731608096609877$ \\
      &$a_2=0.225287493267702165807274831864$ &$b_2=0.255156544139293944162028807345$ \\
      &$a_3=0.416276189612257117795363856737$ &$b_3=0.0181334688208317251361460684041$ \\
      &$a_4=-0.384567270213950399652168569029$ &$b_4=-0.179040110299264554587007062749$ \\
      &$a_5=0.0997271783470514816674547589369$ &$b_5=-0.118470801433302245053382954342$ \\
      &$a_6=-0.108833834399100218757003157958$ &$b_6=0.186461689273821083344937258279$ \\
      &$a_7=0.222010736648991680848341975522$ &$b_7=0.459041581767136840219244627361$ \\
      &$a_8=0.523879522036734296002247438223$ &$b_8=-0.003660836270318358975321459399$ \\
      &$a_9=\frac{1}{2}-\sum_{i=1}^8a_i$ &$b_9=1-2\sum_{i=1}^8b_i$ \\ 
      \cline{2-3}
      &\\
      &\multicolumn{1}{c}{$a_i$}&\multicolumn{1}{c}{$b_i$} \\ 
      \cline{2-3}
      $\mathcal{A}_{18}$&$a_1=0.0866003822712445920135805954462$ &$b_1=-0.08$ \\
      &$a_2=-0.0231572735424388070228714693753$ &$b_2=0.209460550048243262121199483001$ \\
      &$a_3=0.191410576083774088999564416369$ &$b_3=0.274887805875735483503233064415$ \\
      &$a_4=0.378895558692931579545387584925$ &$b_4=-0.224214208870409561366168655624$ \\
      &$a_5=-0.0467359566364556111599485526051$ &$b_5=0.347657740563761656321390026010$ \\
      &$a_6=-0.156198111997810415438979605642$ &$b_6=-0.168783183866211679175007668385$ \\
      &$a_7=0.156025836895094823718831871041$ &$b_7=0.144209344805460873709120777707$ \\
      &$a_8=0.252844012473796333586850465807$ &$b_8=0.0116851121360265483381405054244$ \\
      &$a_9=-0.640644212172254239866860564270$ &$b_9=\frac{1}{2}-\sum_{i=1}^8b_i$ \\
      &$a_{10}=1-2\sum_{i=1}^9a_i$ & \\
      \cline{2-3}
      &\\
      &\multicolumn{1}{c}{$a_i$}&\multicolumn{1}{c}{$b_i$} \\ 
      \cline{2-3}
      $\mathcal{A}_{19}$&$a_1=0.0505805$ &$b_1=0.129478606560536730662493794395$ \\
      &$a_2=0.149999$ &$b_2=0.222257260092671143423043559581$ \\
      &$a_3=-0.0551795510771615573511026950361$ &$b_3=-0.0577514893325147204757023246320$ \\
      &$a_4=0.423755898835337951482264998051$ &$b_4=-0.0578312262103924910221345032763$ \\
      &$a_5=-0.213495353584659048059672194633$ &$b_5=0.103087297437175356747933252265$ \\
      &$a_6=-0.0680769774574032619111630736274$ &$b_6=-0.140819612554090768205554103887$ \\
      &$a_7=0.227917056974013435948887201671$ &$b_7=0.0234462603492826276699713718626$ \\
      &$a_8=-0.235373619381058906524740047732$ &$b_8=0.134854517356684096617882205068$ \\
      &$a_9=0.387413869179878047816794031058$ &$b_9=0.0287973821073779306345172160211$ \\
      &$a_{10}=\frac{1}{2}-\sum_{i=1}^9a_i$ &$b_{10}=1-2\sum_{i=1}^9b_i$ \\
      \cline{2-3}
    \end{tabular}
    \caption{\small Coefficients of 8th-order RKN splitting methods of type $\mathcal{A}_{s}$, with $s=17$, 18 and 19 stages.}
    \label{tau:rkna}
  \end{center}
  }
\end{table}

\paragraph{$s=18$ stages.} 
With one more stage we have one free parameter that can be used to get in principle smaller values of the effective error and eventually more efficient schemes,
as is common in the literature. Notice that the problem in this case involves solving a system of 18 polynomial equations
with 19 variables. Our strategy is the following: for a composition of type $\mathcal{A}_s$ with $s=18$,
we take $a_1$ as the free parameter, and explore the interval $a_1 \in [0,1]$ (since we are interested
in small values of the coefficients) by fixing each time the value of $a_1$. Starting with $2 \times 10^6$ initial points, we have found 722 valid
solutions, the most promising corresponding to the choice $a_1 = 0.08$. This solution is then taken as the starting point of an arc-length continuation method
and follow the solution along the curve leading to a local minimum of the 1-norm of the vector of coefficients. In doing so we apply the algorithm presented
in \cite{makazaga00ac,alberdi19aab}.
After this process, we check several methods
in practice and finally the solution $\mathcal{A}_{18}$ collected in Table \ref{tau:rkna}, with  $E_{\mathrm{f}}$, $\Delta$ and $\delta$ 
given in Table \ref{tau:ef}.

The same technique is applied to compositions $\mathcal{B}_{18}$ leading to the solution collected in Table 
\ref{tau:rknb} after 1070748 initial points and the application of arc-length continuation.

\paragraph{$s=19$ stages.} Adding an additional stage and so forming the composition $\mathcal{A}_{19}$, we have explored the space of parameters 
in the region $a_1, a_2 \in [0.05, 0.15]$, where we have found 295 valid solutions. Then, we start from the one with best parameters and apply the
following strategy: let us denote by $\mathbf{u}_0$ the vector of coefficients of this initial solution. Then we generate a random vector $\mathbf{\alpha}$
verifying $\mathbf{\alpha} \cdot (\mathbf{u} - \mathbf{u}_0) = 0$. Now we apply continuation along the curve that results from the intersection of the space
of solutions (with 2 free parameters) with the random generated hyperplane. The final solution is collected in Table \ref{tau:rknb}.

Concerning the composition $\mathcal{B}_{19}$, 173 solutions have been obtained out of more than $1.3 \times 10^6$ initial points. After applying the
previous technique, we arrive at the solution reported in Table \ref{tau:rknb}.

\begin{table}[!h]
{\small
  \renewcommand\arraystretch{1.1}
  \begin{center}
    \begin{tabular}{lll}
      &\multicolumn{1}{c}{$a_i$}&\multicolumn{1}{c}{$b_i$} \\ 
      \cline{2-3}
      $\mathcal{B}_{17}$&$a_1=0.160227696073839513690970240076 $ &$b_1=0.0514196142537210073343152693459$ \\
      &$a_2=0.306354507436867319879440957100$ &$b_2=0.250497030318342871458417941091$ \\
      &$a_3=0.308395508895171191756544975556$ &$b_3=0.512412268300327350035492806653$ \\
      &$a_4=0.120362086566233408450063177659$ &$b_4=-0.231597138650894401279645184364$ \\
      &$a_5=-0.622888687549183872072186218718$ &$b_5=0.116091323536875759881216298975$ \\
      &$a_6=0.635560951632990078378672016548$ &$b_6=-0.0098365173246965763985763034283$ \\
      &$a_7=-0.144226974795419229640437363913$ &$b_7=-0.108032771466281638634277563747$ \\
      &$a_8=-0.284867527074173816678992817545$ &$b_8=0.249039864198023642002940910070$ \\
      &$a_9=1-2\sum_{i=1}^8a_i$ &$b_9=\frac{1}{2}-\sum_{i=1}^8b_i$ \\
      \cline{2-3}
      &\\
      &\multicolumn{1}{c}{$a_i$}&\multicolumn{1}{c}{$b_i$} \\ 
      \cline{2-3}
      $\mathcal{B}_{18}$&$a_1=0.144410089394373457971755553148$ &$b_1=0.045$ \\
      &$a_2=0.911935520865154315536815857376$ &$b_2=0.459016679491512416807266107555$ \\
      &$a_3=-0.00072932909837392655161199996844$ &$b_3=-0.0456553445594333153223655352757$ \\
      &$a_4=-0.930317101800698721159455541447$ &$b_4=0.0457031020401841003192648096559$ \\
      &$a_5=0.253804074671714046593439154323$ &$b_5=-0.216814341025322492810152535338$ \\
      &$a_6=0.147948981530918626913598733391$ &$b_6=0.163168264552484857133047358600$ \\
      &$a_7=-0.448814759614614928125216243784$ &$b_7=-0.0857080319814376219389850039430$ \\
      &$a_8=0.0824123980794580106751237195418$ &$b_8=0.0265745810650523466142922093591$ \\
      &$a_9=\frac{1}{2}-\sum_{i=1}^8 a_i$ &$b_9=-0.0365538332992893220147096150675$ \\
      &  & $b_{10}=1-2\sum_{i=1}^9b_i$  \\
      \cline{2-3}
      &\\
      &\multicolumn{1}{c}{$a_i$}&\multicolumn{1}{c}{$b_i$} \\ 
      \cline{2-3}
      $\mathcal{B}_{19}$&$a_1=0.337548675291317241942440116575$ &$b_1=0.036132460472136313416730168194$ \\
      &$a_2=-0.223647977575409990331768222380$ &$b_2=0.012697863961074113381675193011$ \\
      &$a_3=0.168949714872223740906385138015$ &$b_3=0.201318391240629276109068041836$ \\
      &$a_4=0.171179938816205886154783136334$ &$b_4=0.135683350134504233201330671671$ \\
      &$a_5=-0.349765168067292877221144631312$ &$b_5=-0.0579071833999963041504740663015$ \\
      &$a_6=0.523808861006312397712070357524$ &$b_6=-0.0772509501792649549463874931821$ \\
      &$a_7=-0.194208871063049124066394765282$ &$b_7=-0.00264758266409925952822161203471$ \\
      &$a_8=-0.323496751337931087309823477561$ &$b_8=-0.0329844384945603065320797537355$ \\
      &$a_9=0.322817287614899749216601693799$ &$b_9=0.0476781560950366927530646289755$ \\
      &$a_{10}=1-2\sum_{i=1}^9 a_i $ &$b_{10}=\frac{1}{2}-\sum_{i=1}^9 b_i$ \\
      \cline{2-3}
    \end{tabular}
    \caption{\small Coefficients of 8th-order RKN splitting methods of type $\mathcal{B}_{s}$, with $s=17$, 18 and 19 stages.}
    \label{tau:rknb}
  \end{center}
}  
\end{table}

Although the quantities (\ref{eferror}) and (\ref{sizeco}) provide useful information about the quality and relative performance of the methods,
one should have in mind that the size of the error terms and therefore the efficiency of each scheme ultimately depends on the particular problem one
is considering and even on the initial conditions. 
For this reason it is convenient to check the behavior of the different schemes on a variety of differential equations and initial conditions, and also to compare
them with other efficient numerical integrators available in the literature. We have separated the numerical illustrations into two sections. Thus, in section
\ref{num.test1} we compare the new schemes with symmetric compositions (\ref{eq.2.1.2}) of order 8, whereas in section \ref{num.test2} we also consider 
RKN splitting integrators of orders 4 and 6, as well as extrapolation methods.

\section{Numerical tests I: 8th-order schemes}
\label{num.test1}

The first set of examples is intended to illustrate the performance of the new RKN splitting methods in comparison with the most 
efficient symmetric compositions of the form (\ref{eq.2.1.2}) we have found in the literature. 
In addition, we also include in the tests the only 8th-order RKN splitting method with 17 stages. 
Specifically, in addition to the previous $\mathcal{A}_s$ and $\mathcal{B}_s$ schemes, we consider the following 8th-order integrators:
\begin{itemize}
\item $\mathcal{O}_{17}$: the RKN splitting method of type $\mathcal{A}_s$ presented in \cite{okunbor94oeo}, with $s=17$ stages.
\item $\mathcal{SS}_{17}$: the symmetric composition of $m=17$ symmetric 2nd-order methods of the form (\ref{eq.2.1.2}) obtained in \cite{kahan97ccf} 
(the coefficients are also collected in \cite[p. 157]{hairer06gni}).
\item $\mathcal{SS}_{19}$ and $\mathcal{SS}_{21}$:  schemes (\ref{eq.2.1.2}) with $m=19$  and $m=21$, respectively, presented in
\cite{sofroniou05dos}.
\end{itemize}
These $\mathcal{SS}_{m}$ methods have been shown to be the most efficient 8th-order schemes within the family of composition methods (\ref{eq.2.1.2}).
We collect in Table \ref{tab.ef.2} the corresponding values of the quantities $E_{\mathrm{f}}$ and $\Delta$ for methods
$\mathcal{SS}_{m}$ when they are used with $\mathcal{S}_h^{[2]}$ as in (\ref{strang.a}) (ABA) or (\ref{strang.b}) (BAB).  The values
of $E_{\mathrm{f}}$ are always greater when the basic scheme is (\ref{strang.b}).

The implementation of all the integrators has been done in Python 3.7~\cite{python3rm} running on Debian GNU/Linux 10~\cite{krafft2005debian}
and the array operations have been coded using the \textit{NumPy} library~\cite{numpy}.

\begin{table}[!h]
  \begin{center}
    \begin{tabular}{lllllll}
      &&\multicolumn{3}{c}{$E_{\rm{f}}$}  &&\multicolumn{1}{c}{$\Delta$} \\
      \hline
       &&ABA &&BAB & \\ \hline
      $\mathcal{O}_{17}$ &\hspace{2.5cm} &4.78 && -- &\hspace{1cm} &16.63 \\ 
      $\mathcal{SS}_{17}$ &&3.12 &&3.30 &&8.33 \\ 
      $\mathcal{SS}_{19}$ &&2.66 &&2.68 &&6.84 \\ 
      $\mathcal{SS}_{21}$ &&2.59 &&2.88 &&6.43 \\ 
      \hline
    \end{tabular}
    \caption{{\small Effective error  $E_{\mathrm{f}}$ and 1-norm of the vector of coefficients for 8th-order symmetric compositions of symmetric methods
    $\mathcal{SS}_{m}$ and the RKN splitting method of \cite{okunbor94oeo}.}}
    \label{tab.ef.2}
  \end{center}
\end{table}

\begin{figure}[!h]
  \begin{center}
    \subfloat[]{
      \includegraphics[width=8.2cm]{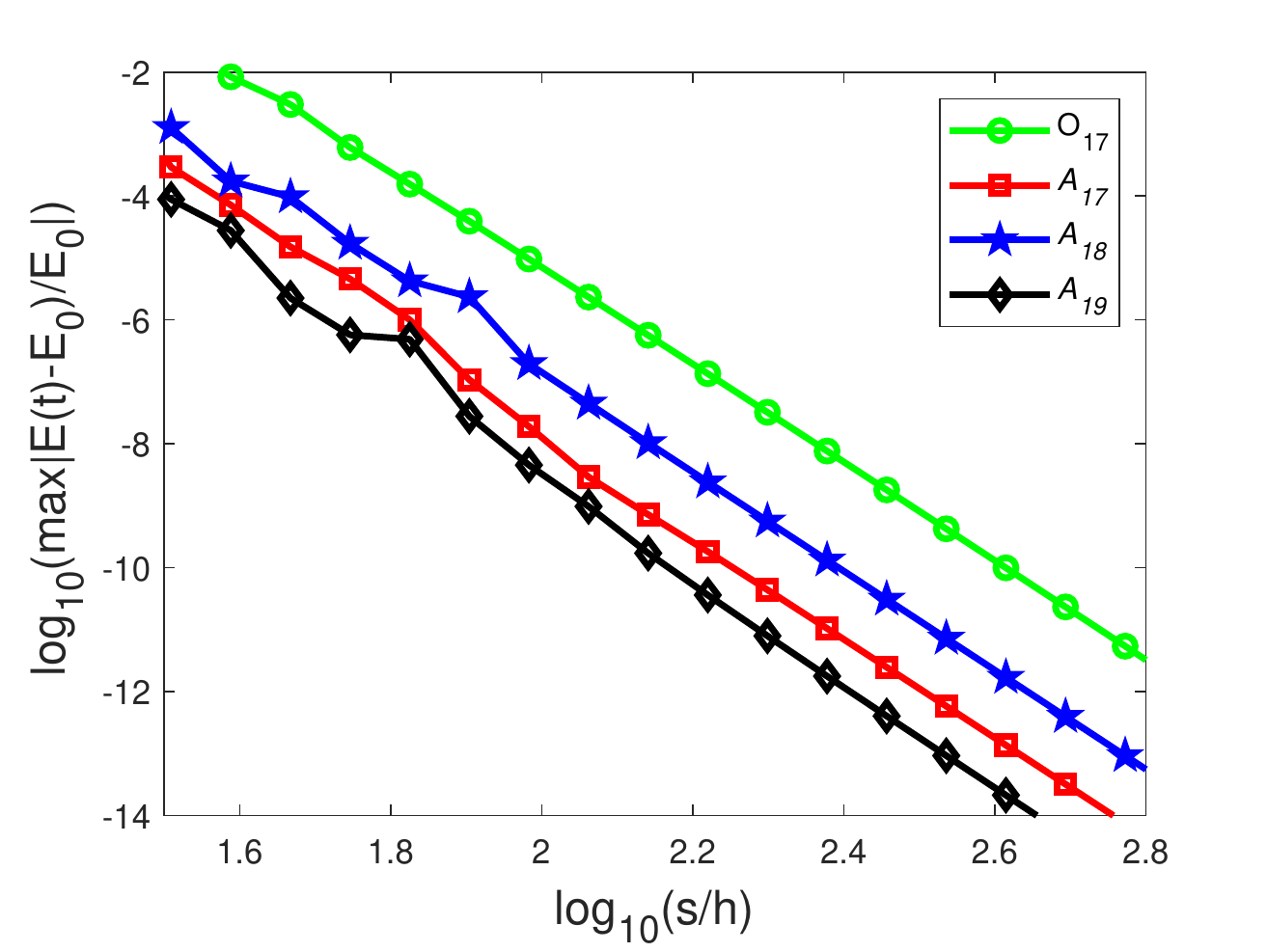}}
    \subfloat[]{
      \includegraphics[width=8.3cm]{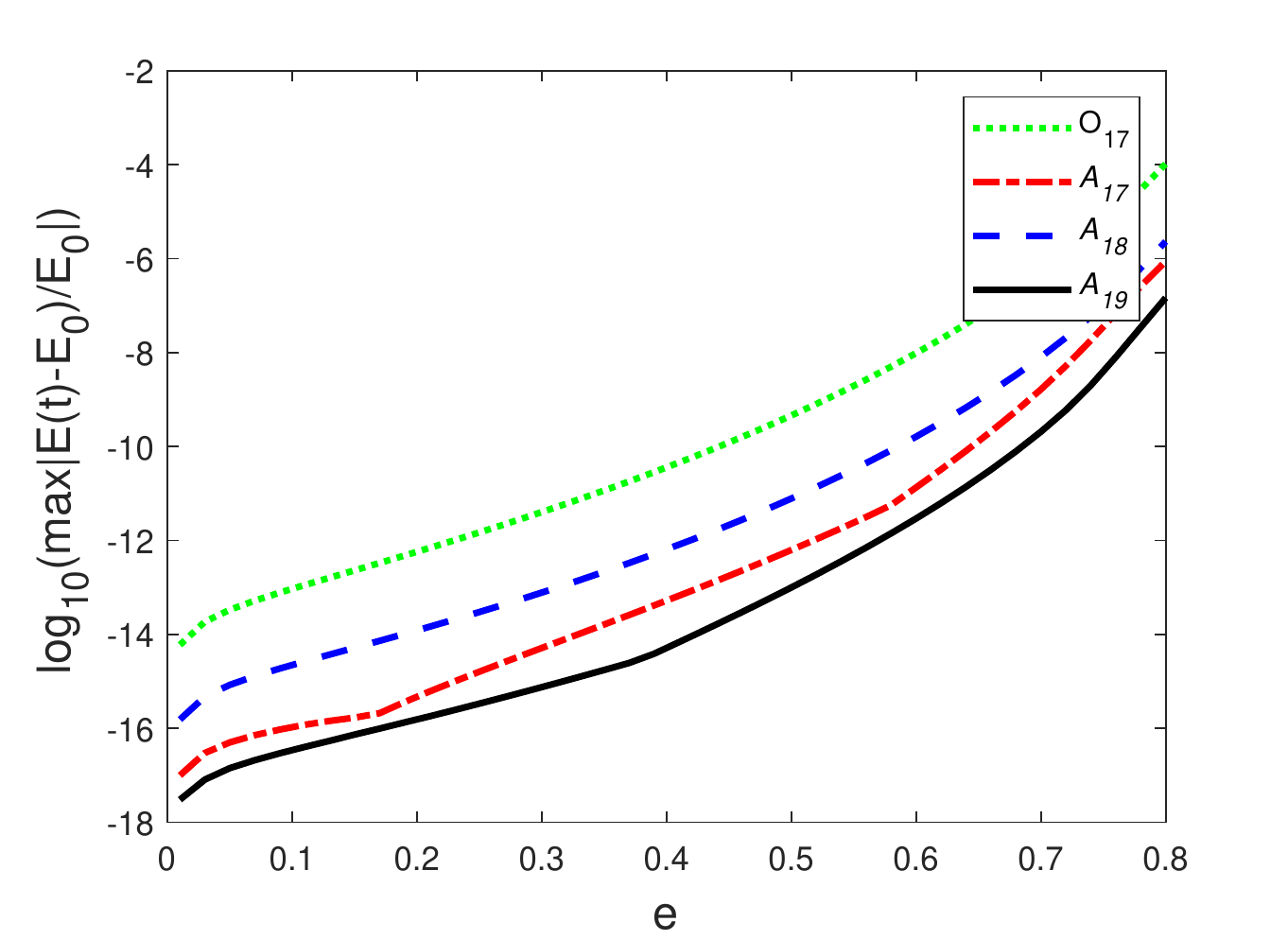}}
    \caption{\small (a) Efficiency diagram for the Kepler problem with $e=0.5$ for all RKN splitting methods of $\mathcal{A}_s$ type. The final time is $t_f=1000$. (b) Maximum error in energy for different values of the eccentricity with $t_f=1000$ and $s/h=340$.}
    \label{fig:kepler_rkn}
  \end{center}
\end{figure}

\begin{figure}[!h]
  \begin{center}
    \subfloat[]{
      \includegraphics[width=8.5cm]{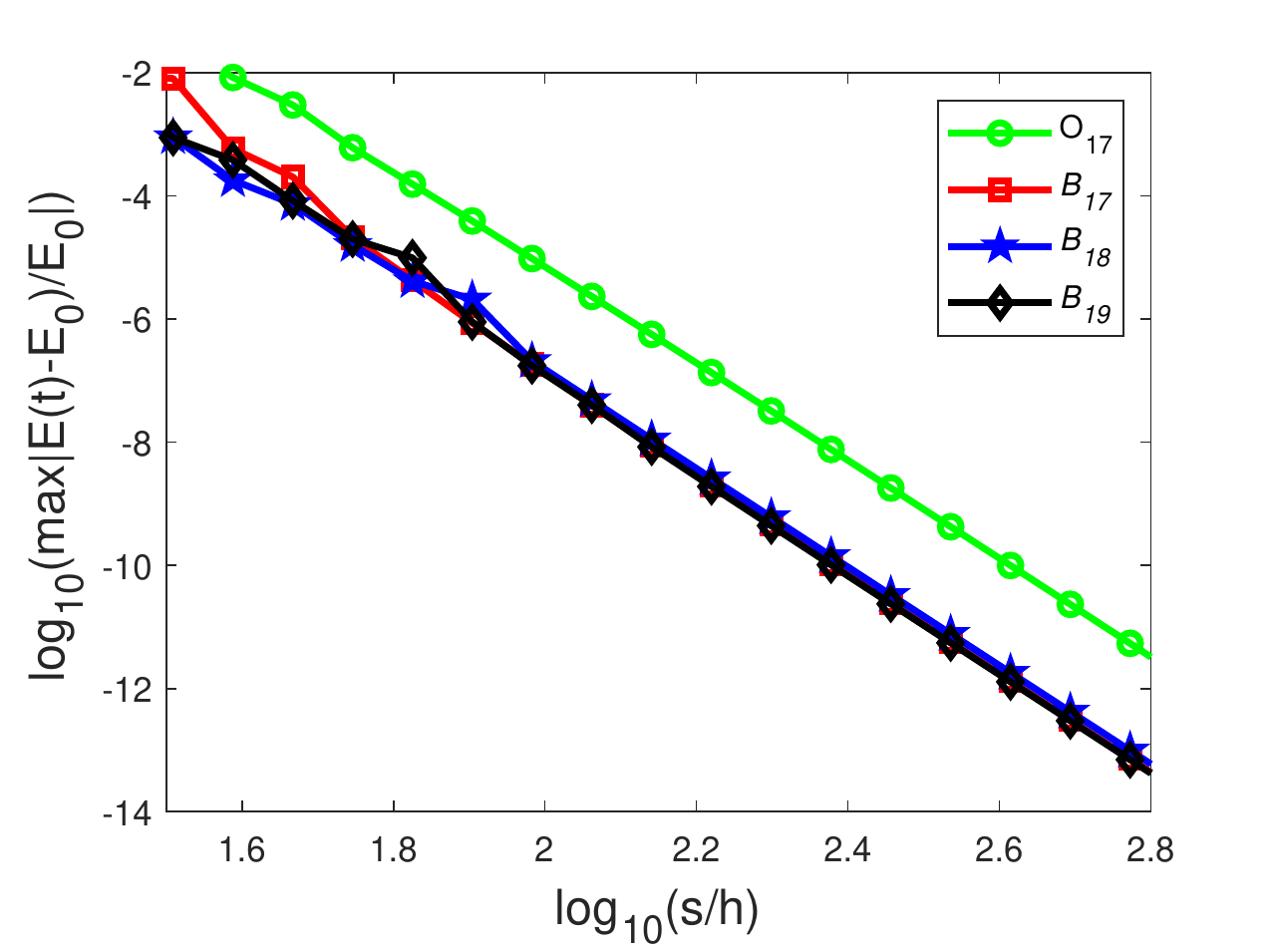}}
    \subfloat[]{
      \includegraphics[width=8.1cm]{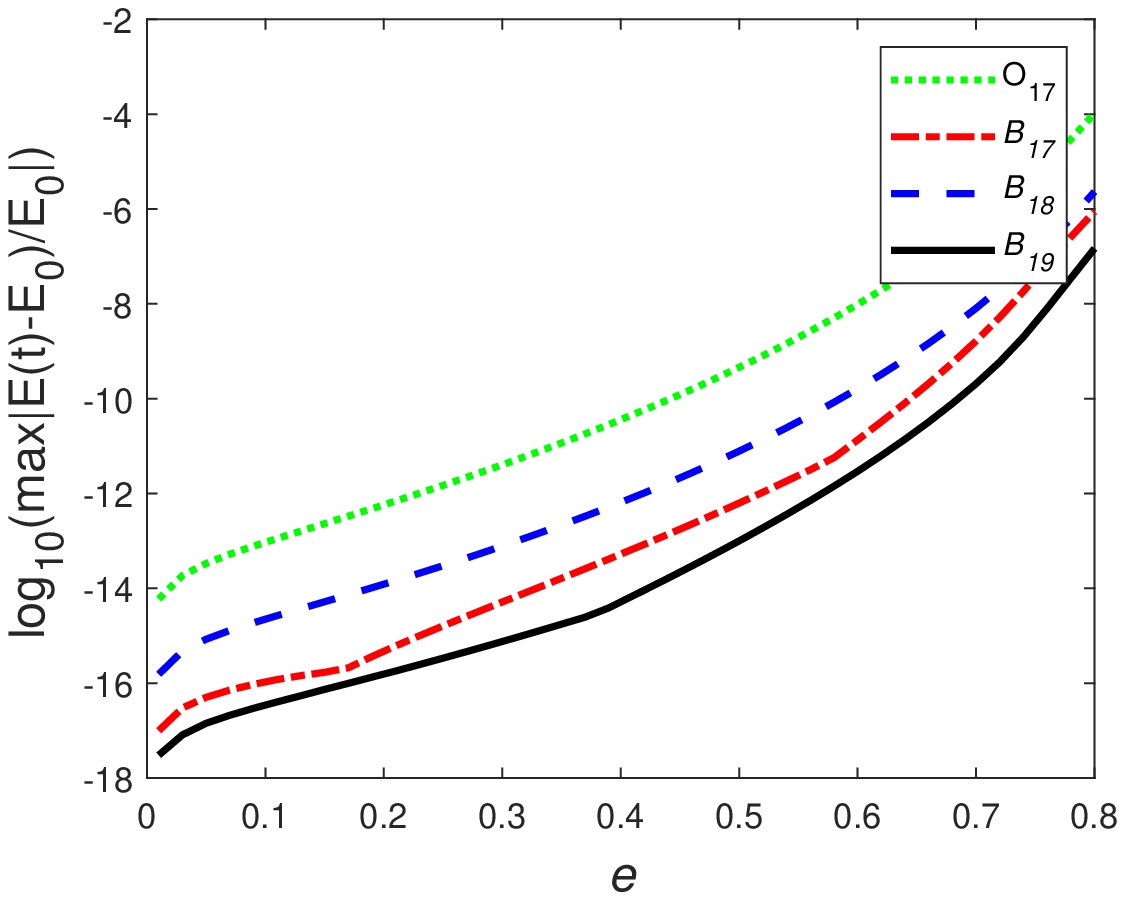}}
    \caption{\small (a) Efficiency diagram for the Kepler problem with $e=0.5$ for all RKN splitting methods of $\mathcal{B}_s$ type. The final time is $t_f=1000$. (b) Maximum error in energy for different values of the eccentricity with $t_f=1000$ and $s/h=340$.}
    \label{fig:kepler_rknb}
  \end{center}
\end{figure}

\paragraph{Example 1: Kepler problem.}
We take the 2-body gravitational problem with Hamiltonian
\begin{equation} \label{eq.HamKepler}
   H(q,p) = 
	\frac{1}{2} p^T p - \mu \frac{1}{r},
\end{equation}
where $q=(q_1,q_2), p=(p_1,p_2)$, $\mu=GM$, $G$ is the gravitational constant and $M$ is the sum of the masses of the two bodies. 
We take $\mu=1$ and initial conditions
\begin{equation} \label{eq.InitKepler}
  q_1(0) = 1- e, \quad q_2(0) = 0, \quad p_1(0) = 0, \quad p_2(0) = \sqrt{\frac{1+e}{1-e}},
\end{equation}
so that the trajectory corresponds to an ellipse of eccentricity $e$, with period $2 \pi$ and energy $E=-\frac{1}{2}$. We first check the order of the
new RKN splitting methods and compare their efficiency with respect to $\mathcal{O}_{17}$. Thus, Figure \ref{fig:kepler_rkn} (left panel)
shows the relative error
in energy with respect to  $s/h$ (which is proportional to the number of force evaluations) for $e=0.5$ and a final time $t_f=1000$ for methods of type $\mathcal{A}_s$, 
whereas in the right panel we explore the range of eccentricities $0 \le e \le 0.8$.
All schemes involve the same number of evaluations of the potential in this case. Figure \ref{fig:kepler_rknb} shows analogous results for methods of type
$\mathcal{B}_s$. 
Notice that the order 8 is clearly visible in the figures and that the 
new methods are more efficient than $\mathcal{O}_{17}$. 
The improvement is particularly prominent for $\mathcal{A}_{17}$ and specially $\mathcal{A}_{19}$ 
(up to four orders of magnitude for the same value of $h/s$) and is more moderate for methods $\mathcal{B}_s$. In fact, all of them show
essentially the same performance, which is lower than that of $\mathcal{A}_{19}$. 

\begin{figure}[!h]
  \begin{center}
    \subfloat[]{
      \includegraphics[width=8.cm]{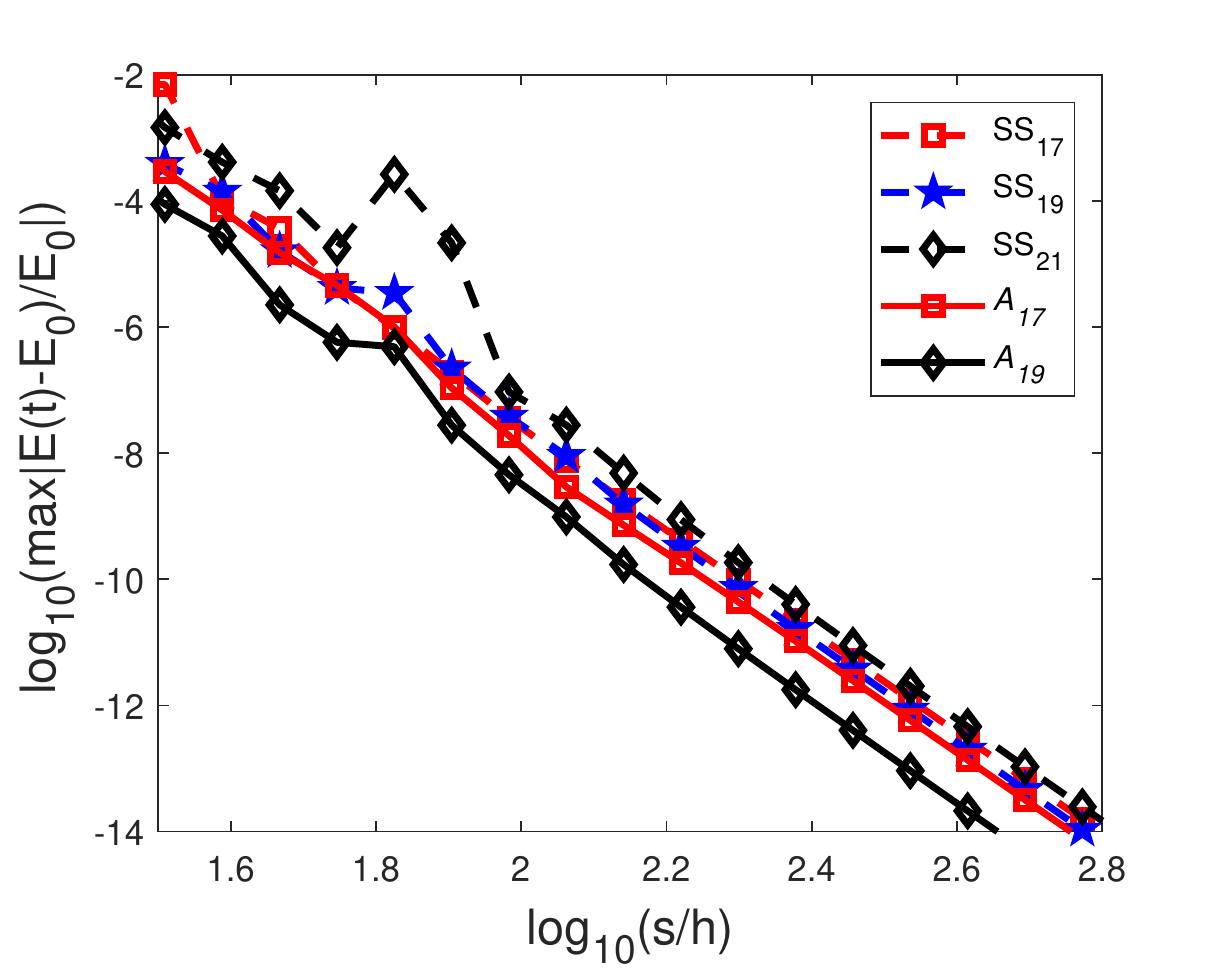}}
    \subfloat[]{
      \includegraphics[width=8.6cm]{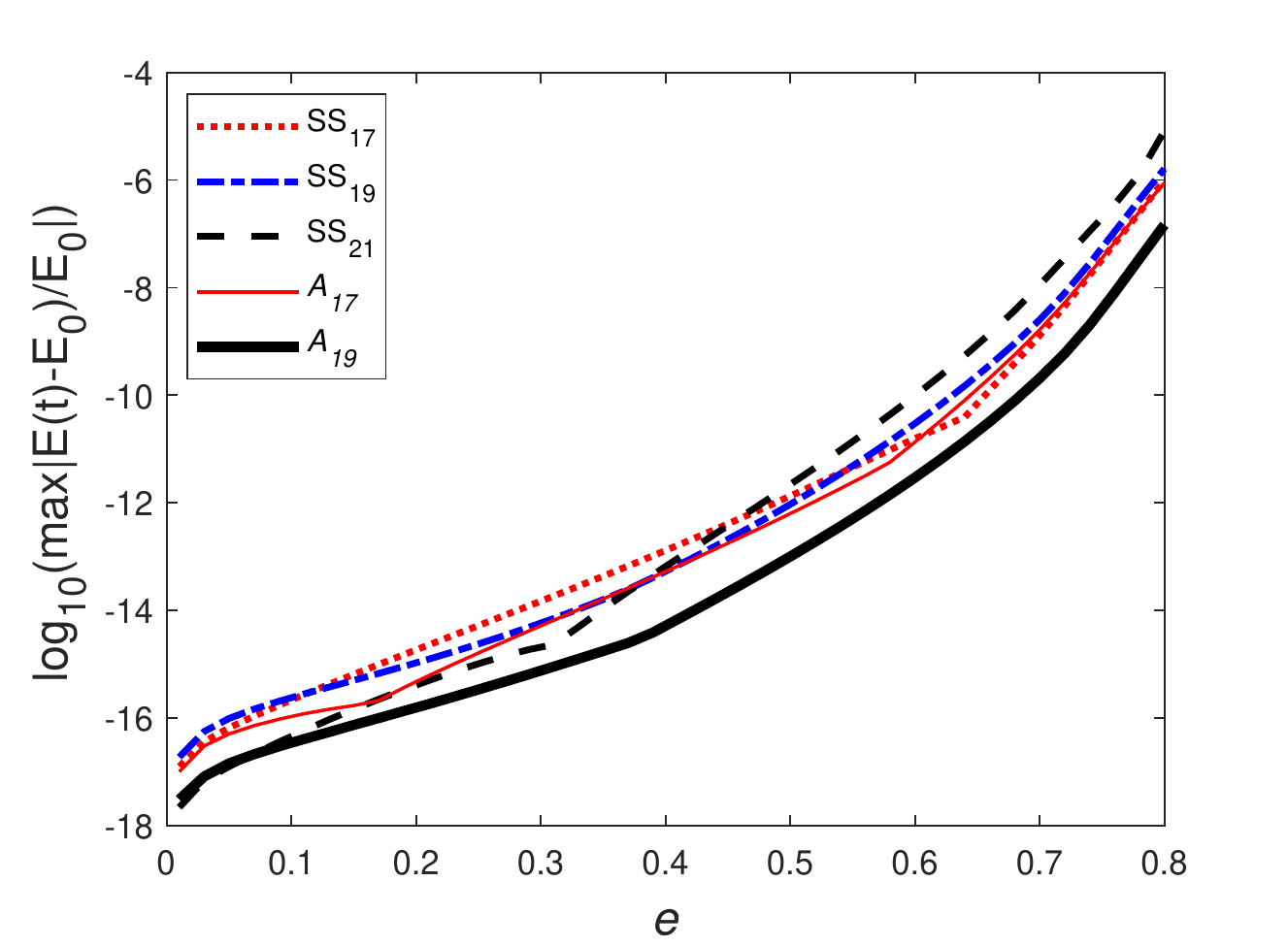}}
    \caption{\small (a) Efficiency diagram for the Kepler problem with $e=0.5$ for composition $\mathcal{SS}_m$ and the new RKN splitting methods 
    $\mathcal{A}_{17}$ and $\mathcal{A}_{19}$. Final time $t_f=1000$. (b) Maximum error in energy for different values of the eccentricity with $t_f=1000$ and 
    $s/h=340$.}
    \label{fig:kepler_ssrkn}
  \end{center}
\end{figure}

We next carry out the same experiment, but in this case we compare the performance of the new schemes $\mathcal{A}_{17}$ and $\mathcal{A}_{19}$ with 
the previous state-of-the-art symmetric compositions of the Strang splitting $\mathcal{SS}_m$, $m=17, 19, 21$. We take the composition \eqref{strang.a} as the basic $\mathcal{S}_h^{[2]}$ method because it shows the best performance in the numerical experiments. The corresponding results are shown in
Figure \ref{fig:kepler_ssrkn}. We notice that $\mathcal{A}_{19}$ is the more efficient method for the whole range of eccentricities explored.

\begin{figure}[!h]
  \begin{center}
    \subfloat[]{
      \includegraphics[width=8.2cm]{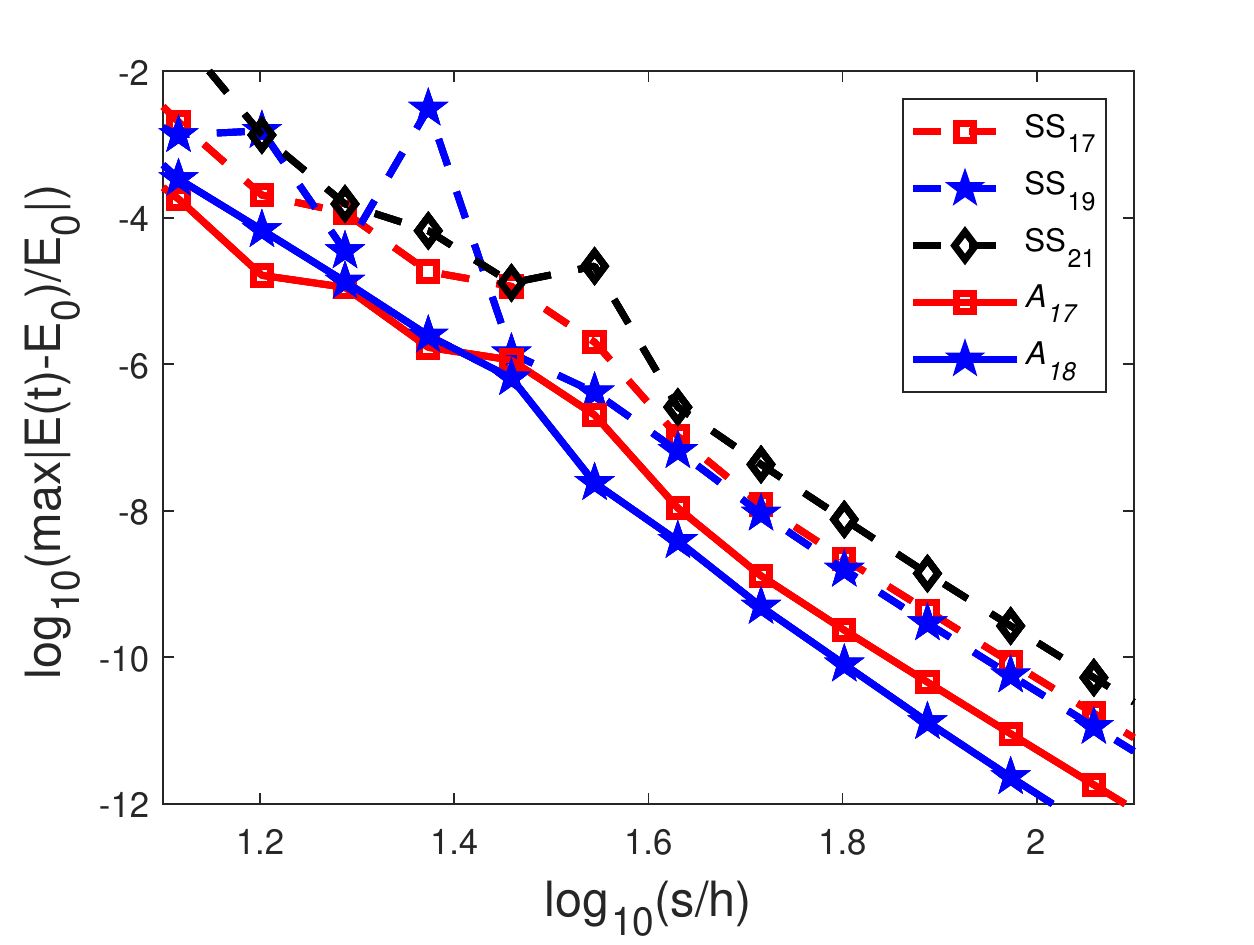}}
    \subfloat[]{
      \includegraphics[width=8.3cm]{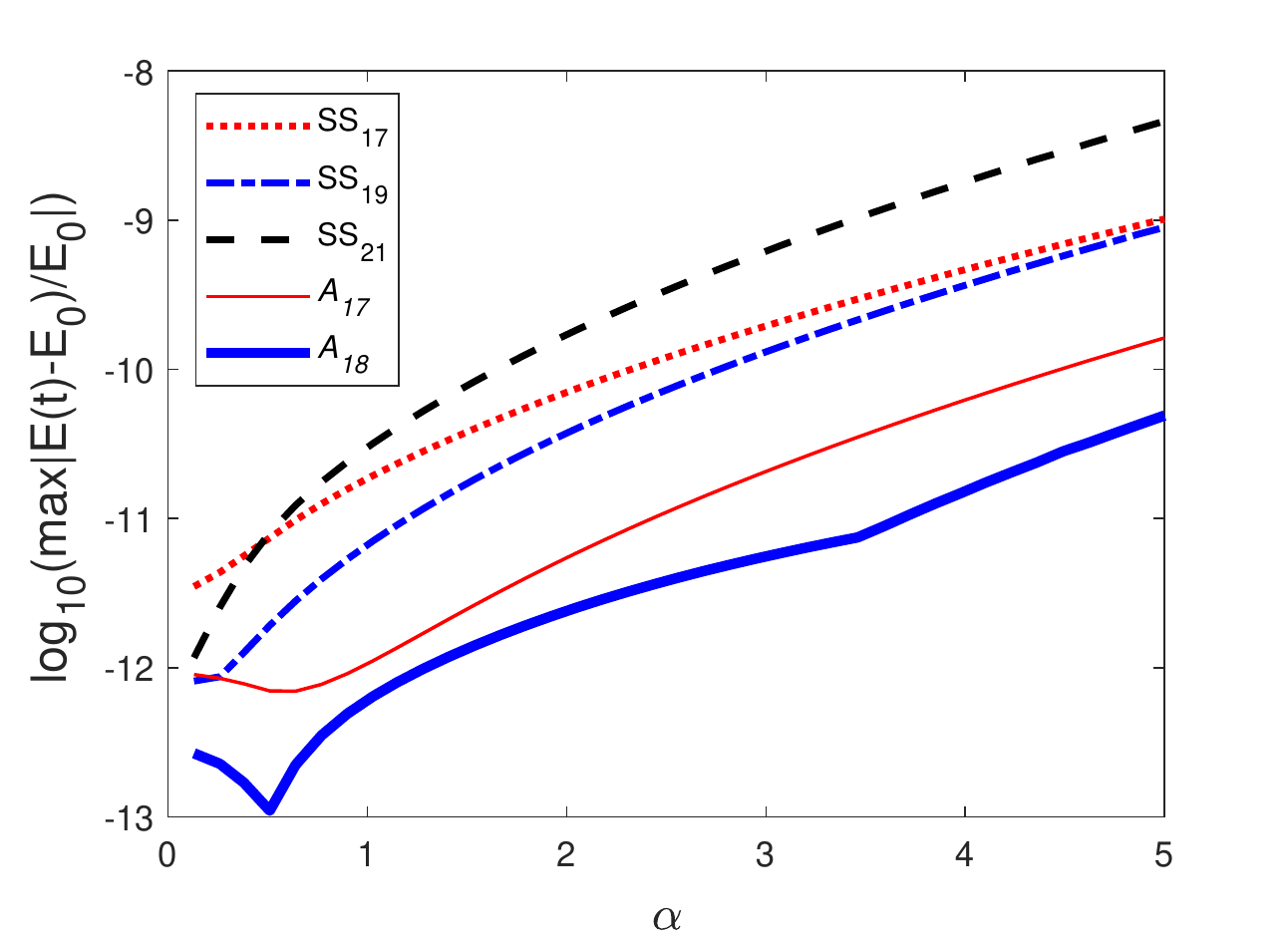}}
    \caption{\small Simple pendulum. (a) Efficiency diagram for $\alpha=3.0$ and final time $t_f=1000$. (b) Maximum error in energy for initial conditions $(q_0,p_0)=(0,\alpha)$ for $\mathcal{SS}$ and the best RKN methods at final time $t_f=1000$ with $s/h=85$.}
    \label{fig:pendol_ssrkn}
  \end{center}
\end{figure}

\paragraph{Example 2: simple pendulum.}
Our next example is the simple mathematical pendulum. In appropriate units, it corresponds to the 1-degree-of-freedom Hamiltonian system
with
\begin{equation} \label{pendulum1}
  H(q,p) = \frac{1}{2} p^2 - \cos q.
\end{equation}
We explore the set of initial conditions $(q,p) = (0, \alpha)$, with $0 \le \alpha \le 5$, integrate until the final time $t_f=1000$ and check the error
in energy along the integration. Since the error achieved by $\mathcal{O}_{17}$ is always 3-4 orders of magnitude larger than the new schemes, we
no longer include them in the diagrams, so that we only compare with symmetric compositions $\mathcal{SS}_m$. Figure \ref{fig:pendol_ssrkn}
shows the efficiency diagram corresponding to $\alpha = 3$ (panel (a)) and the maximum of the relative error in the energy along the integration interval.
In this case, the new schemes $\mathcal{A}_{17}$ and $\mathcal{A}_{18}$ are the most efficient. Scheme $\mathcal{A}_{19}$ shows a similar behavior as
$\mathcal{SS}_{19}$, and thus it has not been included in the diagrams. On the other hand, the most efficient scheme of the BAB type in this case is
$\mathcal{B}_{18}$ (not shown), providing similar results as $\mathcal{A}_{18}$.

\paragraph{Example 3: H\'enon--Heiles potential.}

\begin{figure}[!h]
  \begin{center}
    \subfloat[]{
      \includegraphics[width=8.35cm]{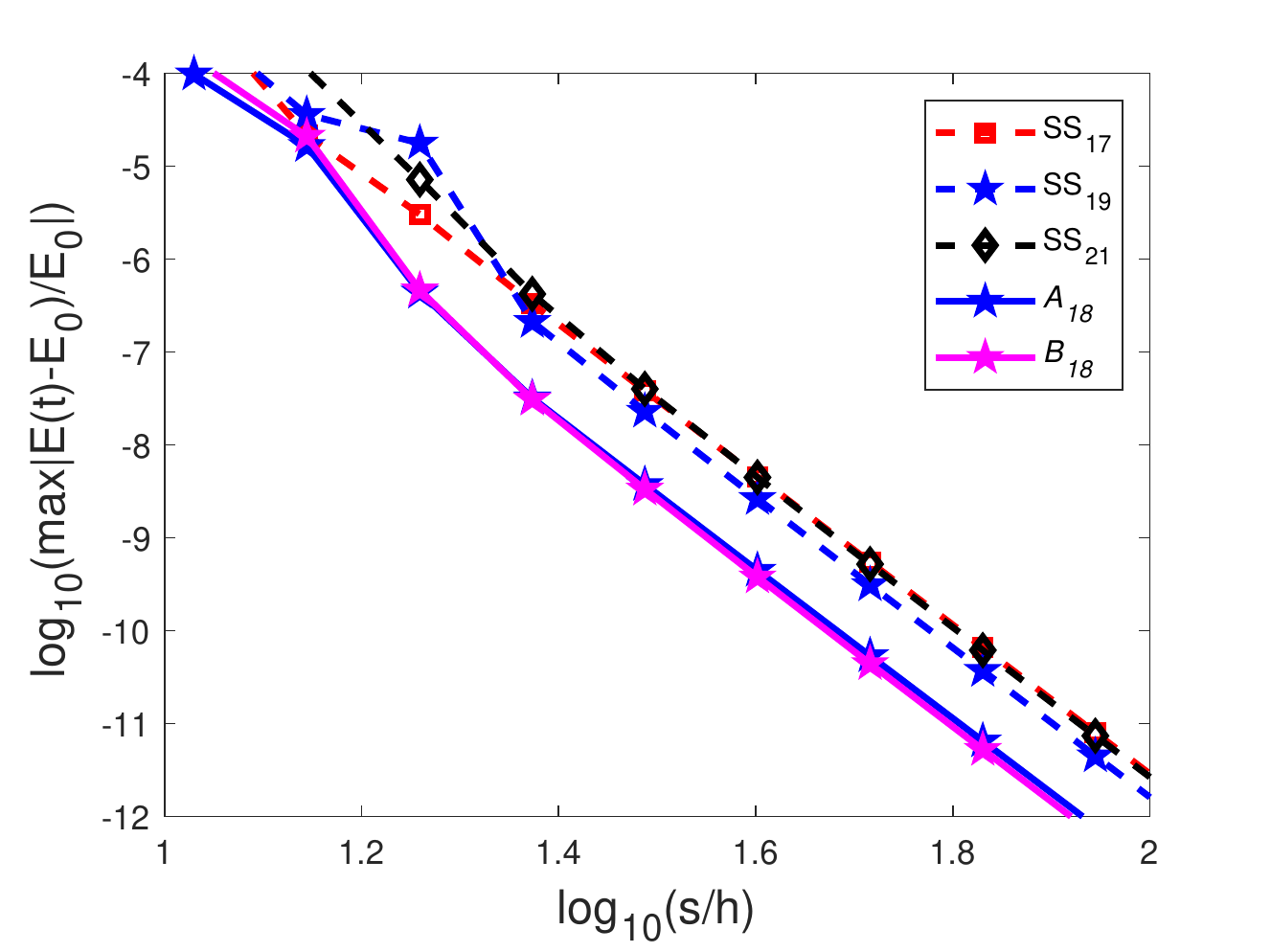}}
    \subfloat[]{
      \includegraphics[width=8.3cm]{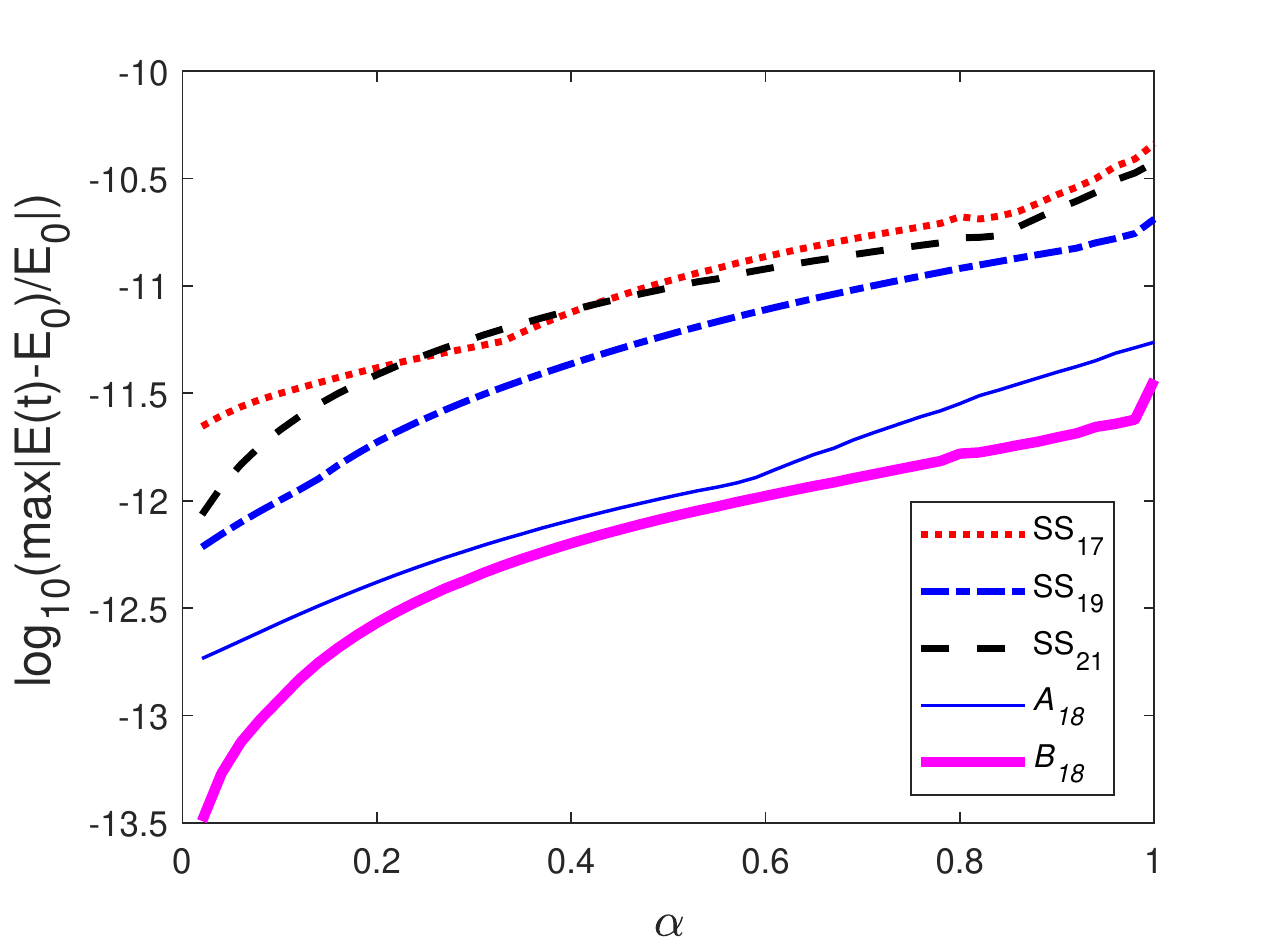}}
    \caption{\small H\'enon--Heiles Hamiltonian. (a) Efficiency diagram with initial condition corresponding to $\alpha=0.2$ and final time $t_f=1000$.
    (b) Maximum error in energy for $0 \le \alpha \le 1$ at final time $t_f=1000$ with $s/h=85$.}
    \label{fig:hehe_ssrkn}
  \end{center}
\end{figure}

For our next experiment we choose the well known two-degrees of freedom H\'enon--Heiles Hamiltonian \cite{henon64tao}
\begin{equation} \label{hh.1}
  H = \frac{1}{2} (p_1^2 + p_2^2) + \frac{1}{2} (q_1^2 + q_2^2) + q_1^2 q_2 - \frac{1}{3} q_2^3.
\end{equation}
It has been the subject of extensive numerical experimentation and 
is considered, in particular, as a model problem to characterize the transition to Hamiltonian chaos. In this case we take the same 
initial conditions as in \cite{blanes02psp}, the set $(q_1,q_2,p_1,p_2) = (\alpha/2, 0, 0, \alpha/4)$, with $0 \le \alpha \le 1$. The
corresponding results are depicted in Figure \ref{fig:hehe_ssrkn}. 
In this case $\mathcal{B}_{18}$ and $\mathcal{A}_{18}$ are the most efficient schemes, whereas $\mathcal{A}_{17}$ is similar as $\mathcal{A}_{18}$ and it is not shown in the figure.

\paragraph{Example 4: Schr\"odinger equation with P\"oschl--Teller potential.}

\begin{figure}[!h]
  \begin{center}
    \includegraphics[width=10cm]{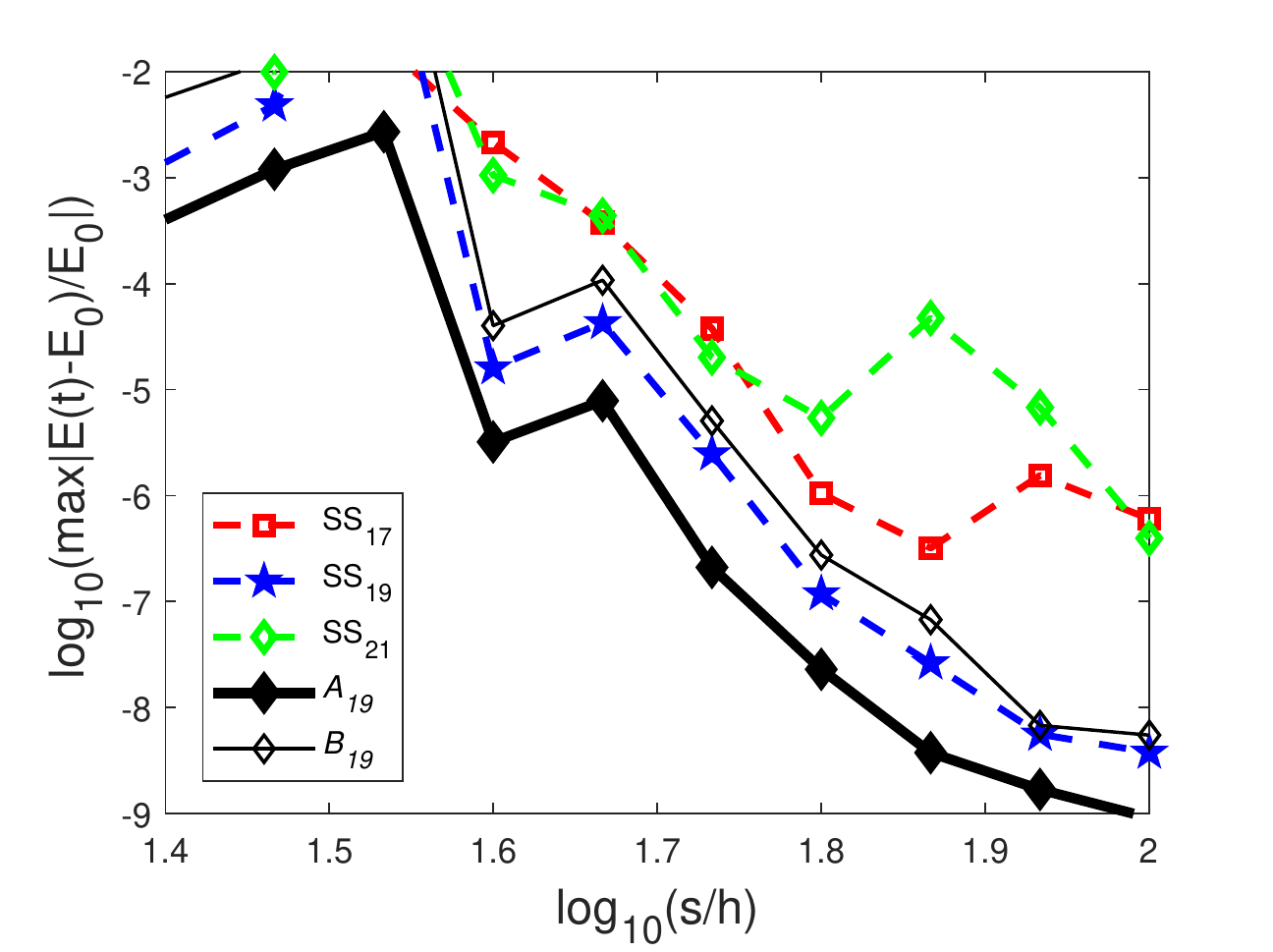}
    \caption{Efficiency diagram of different methods. Schr\"odinger equation with   P\"oschl--Teller potential.}
    \label{fig:pote_ssrkn}
  \end{center}
\end{figure}

Finally, we apply our integrators to the one-dimensional Schr\"odinger equation ($\hbar = 1$)
\begin{equation} \label{schro.1}
  i \frac{\partial}{\partial t} \psi (x,t)  = -\frac{1}{2}
  \frac{\partial^2}{\partial x^2} \psi (x,t) + V(x) \psi (x,t),
\end{equation}
with the well known P\"oschl--Teller potential \cite{flugge71pqm},
\begin{equation} \label{pt1}
  V(x) = -\frac{\lambda (\lambda+1)}{2} \mbox{sech}^2(x),
\end{equation}
with $\lambda(\lambda+1) = 10$. When a Fourier spectral collocation method is used for discretizing in space \cite{trefethen00smi}, one ends up
with the $N$-dimensional linear ODE
\begin{equation}   \label{lin1}
  i \frac{d }{dt} u(t) = H \, u(t)  \equiv (T + V) \, u(t),  \qquad
    u(0) = u_{0} \in \mathbb{C}^N,
\end{equation}
where $T$ is a (full) differentiation matrix related with the kinetic energy, $V$  is a diagonal matrix associated with
the potential and the components of the vector $u$ are the approximations to the wave function at the nodes,
$u_n \approx \psi(x_n,t)$. The action of $T$ on the wave function vector $u$ is then carried out by the forward and backward discrete
Fourier transform (computed with the FFT algorithm) \cite{lubich08fqt}. The initial condition is taken as $\psi_0(x) = \sigma \, \e^{-x^2/2}$, 
with $\sigma$ a normalizing constant,  the interval is 
$x \in [-8,8]$ with $N= 256$ nodes, and the integration is done until the final time $t_f=1000$. In this case we check the error in the expected value of the
energy,
\begin{equation} \label{eq.5.1b}
   \mbox{energy error:} \quad |u_{\mathrm{ap}}^*(t) \cdot (H  u_{\mathrm{ap}}(t))  - u_0^* \cdot (H  u_0)|,
\end{equation}   
where $u_{\mathrm{ap}}(t)$ stands for the numerical approximation obtained by each method. The results are shown in Figure 
\ref{fig:pote_ssrkn}. Observe that the new RKN splitting method $\mathcal{A}_{19}$ is also the most efficient in this setting.

\section{Numerical tests II: RKN splitting and extrapolation methods}
\label{num.test2}

Given the observed improvement of the new 8th-order RKN splitting methods with respect to the symmetric compositions of a basic 2nd-order symmetric
scheme, it seems appropriate to carry out further comparisons with other lower-order RKN splitting methods when medium to high accuracy is desired.
Specifically we consider the following optimized 4th- and 6th-order methods of type $\mathcal{B}_s$ presented in \cite{blanes02psp}:
\begin{itemize}
\item RKN4$_{6}$: order 4 with 6 stages (the scheme SRKN$_6^b$ in  \cite{blanes02psp}). 
\item RKN6$_{11}$: order 6 with 11 stages (the scheme SRKN$_{11}^b$ in  \cite{blanes02psp}). 
\end{itemize}

On the other hand, extrapolation methods constitute one of the most efficient classes of schemes for the numerical integration of the second order differential equation (\ref{rkn.1}) when high accuracy is required \cite{hairer93sod}. Notice, however, that
in contrast with RKN splitting methods, they do not preserve by construction geometric
properties of the exact solution. To carry out our comparisons, we take \eqref{strang.a} as the symmetric second order basic method (which
corresponds to {\it St\"ormer's rule} \cite[eq. (14.32)]{hairer93sod}) and apply the harmonic sequence to construct by extrapolation schemes of 
orders 4, 6 and 8 with only 3, 6 and 10 stages, respectively. For completeness, the resulting methods can be written explicitly  as
\[
  \Psi_{(r=2k)} = \sum_{\ell=1}^k\alpha_{\ell}^{(k)}\prod_{i=1}^{\ell} {\cal S}_{h/\ell}^{[2]}, \qquad\quad k = 2,3,4,
\]
with $\alpha^{(k)}=(\alpha_{1}^{(k)},\ldots,\alpha_{k}^{(k)})$ and
\begin{equation}
	 \alpha^{(2)} =\left(-\frac{1}{3},\frac{4}{3} \right),   \qquad
 \alpha^{(3)} =\left(\frac{1}{24},-\frac{16}{15},\frac{81}{40} \right),   \qquad
 \alpha^{(4)} =\left(-\frac{1}{360},\frac{16}{45},-\frac{729}{280},\frac{1024}{315} \right).
\end{equation}

\paragraph{Example 5: Kepler problem revisited.}

\begin{figure}[!h]
  \begin{center}
    \subfloat[]{
      \includegraphics[width=8.3cm]{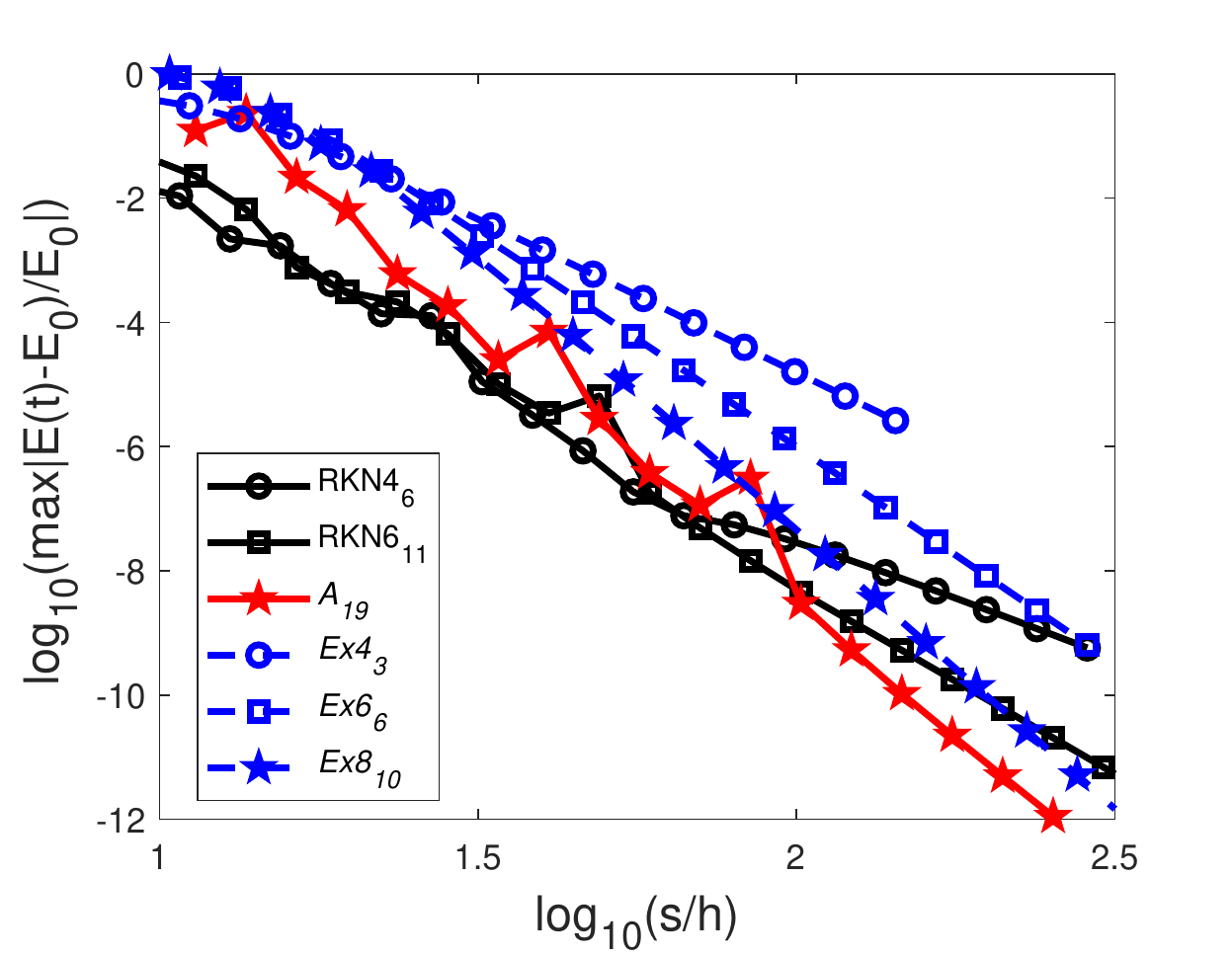}}
    \subfloat[]{
      \includegraphics[width=8.3cm]{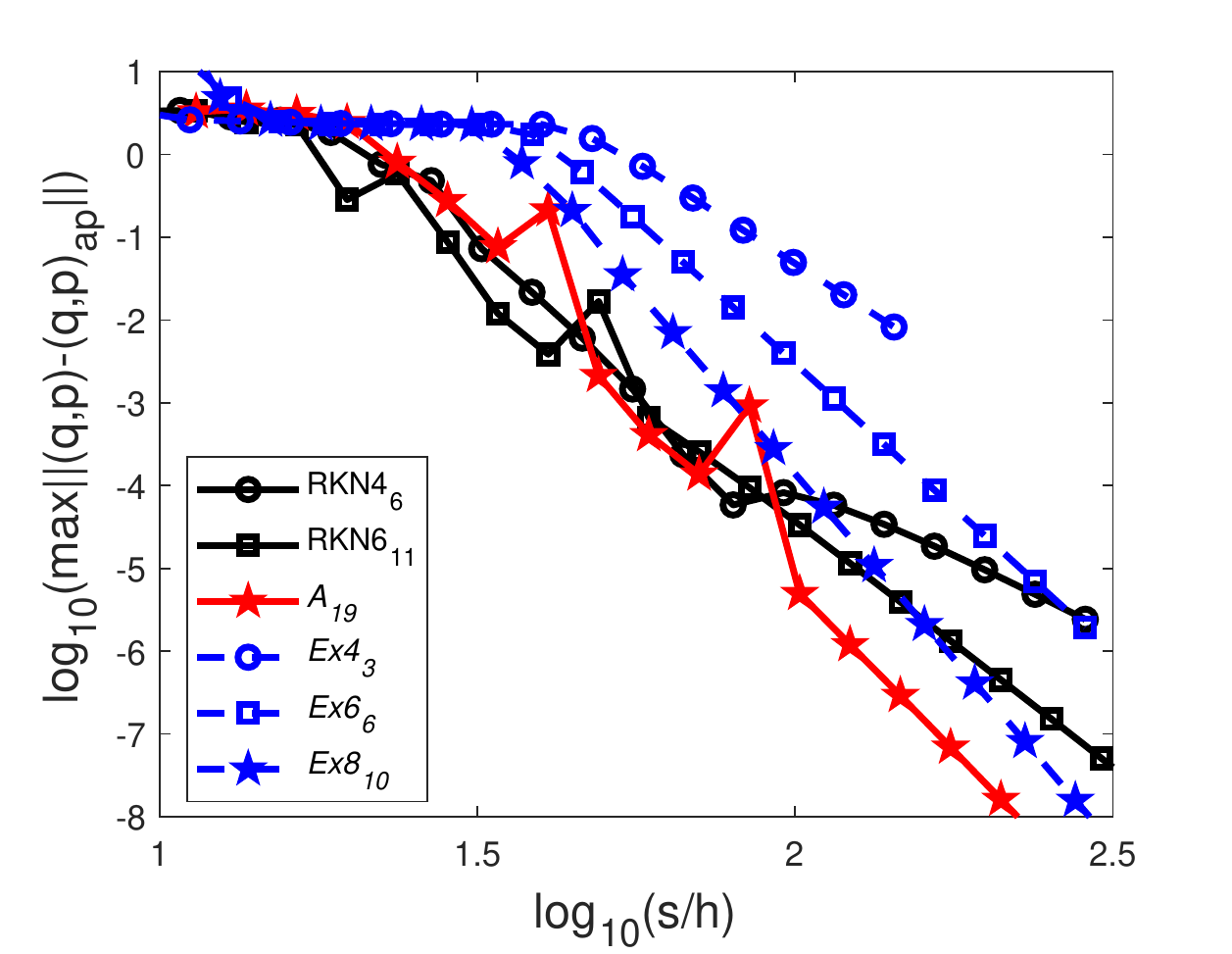}}
    \caption{\small (a) Maximum error in the energy for the Kepler problem with $e=0.5$ obtained by RKN splitting methods RKN4$_{6}$,  RKN6$_{11}$,
    $\mathcal{A}_{19}$ (solid lines), and extrapolation
    (dashed lines) of orders 4 (circles), 6 (squares) and 8 (stars). (b) Same for the maximum error in position.}
    \label{fig:kepler_rknExtrap}
  \end{center}
\end{figure}

For the Hamiltonian (\ref{eq.HamKepler}) with initial conditions (\ref{eq.InitKepler}) we compare the most efficient 8th-order RKN splitting method 
$\mathcal{A}_{19}$ with the 4th- and 6th-order schemes RKN4$_6$ and RKN6$_{11}$, and the previous extrapolation methods of orders 4, 6 and 8 for
the final time $t_f = 1000$.
The results achieved for the maximum error in energy and positions are displayed in Figure \ref{fig:kepler_rknExtrap}. To reduce round-off errors when
computing the linear combinations in extrapolation methods, instead of evaluating directly the numerical solution as $y_{n+1} = \Psi_{(r=2k)} y_n$,
we express $y_{n+1}^{(\ell)} \equiv \prod_{i=1}^{\ell} {\cal S}_{h/\ell}^{[2]} \, y_n$ as $y_{n+1}^{(\ell)} = y_n + \Delta y_{n+1}^{(\ell)}$. In this way 
we compute only $\Delta y_{n+1}^{(\ell)}$, then extrapolation is used only for these increments, namely,
\[
  \Delta y_{n+1} = \sum_{\ell=1}^k\alpha_{\ell}^{(k)}\Delta y_{n+1}^{(\ell)}
\]
and finally we form $y_{n+1}=y_n+\Delta y_{n+1}$. In doing so, round-off errors are reduced by two or more digits.

Figure \ref{fig:kepler_rknExtrap} shows that the new RKN splitting methods are competitive with extrapolation methods and, in particular, $\mathcal{A}_{19}$ is the most efficient when medium to high accuracy is desired.

\begin{figure}[!h]
  \begin{center}
      \includegraphics[width=10cm]{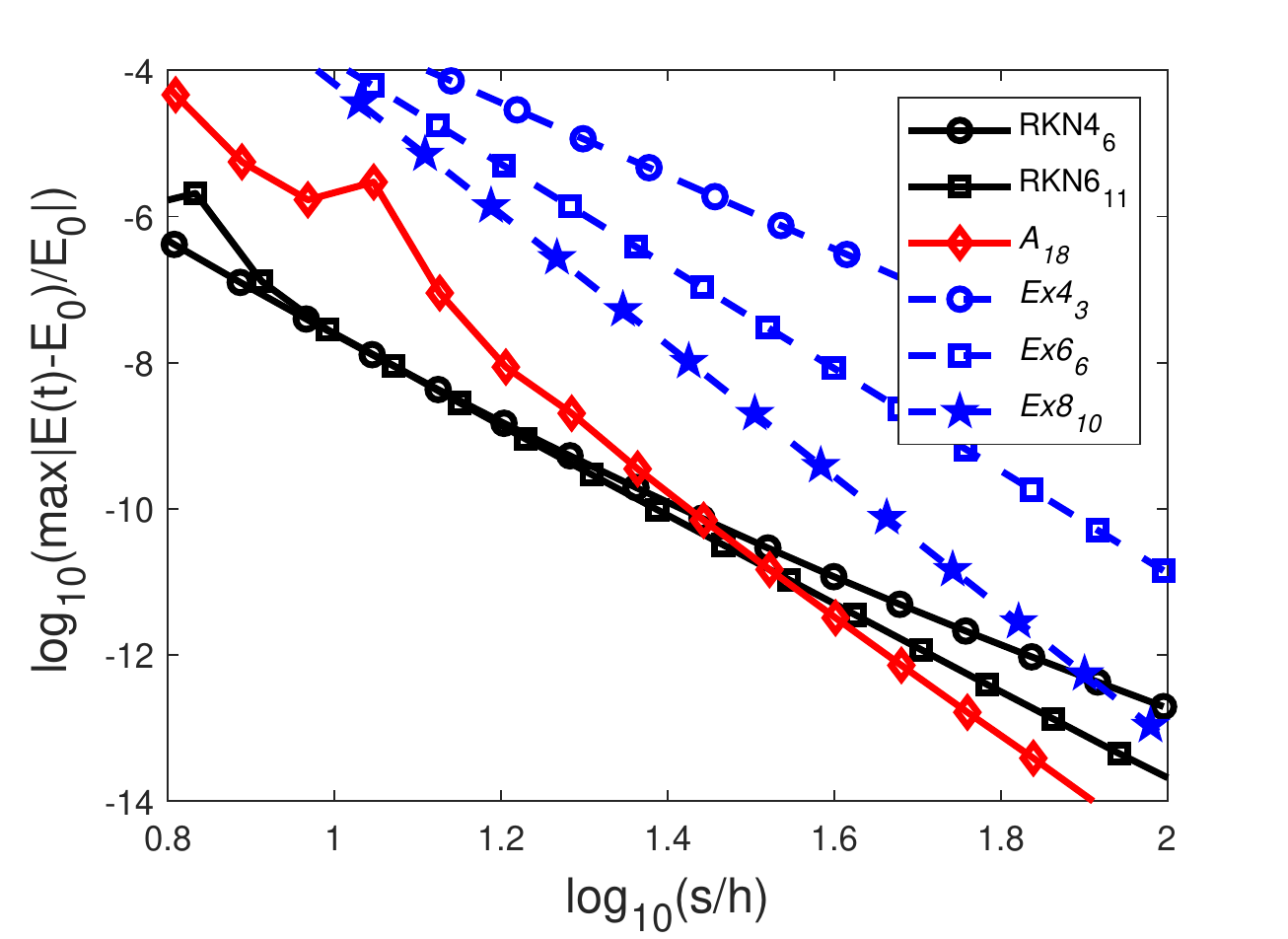}
    \caption{\small Simple pendulum.  Maximum error in the energy for the simple pendulum with initial conditions $(q,p) = (0, 0.3)$ and final time $t_f=1000$ obtained by RKN splitting methods RKN4$_{6}$,  RKN6$_{11}$,
    $\mathcal{A}_{19}$ (solid lines), and extrapolation
    (dashed lines) of orders 4 (circles), 6 (squares) and 8 (stars).}
    \label{fig:pendol_rknExtr}
  \end{center}
\end{figure}

\paragraph{Example 6: simple pendulum revisited.}
Let us consider again the simple pendulum, this time with initial conditions $(q,p) = (0, 0.3)$. We measure the error in energy along the integration for the schemes RKN4$_6$, RKN6$_{11}$, $\mathcal{A}_{18}$ and the extrapolation methods until the final time $t_f=1000$. Figure \ref{fig:pendol_rknExtr}
shows the efficiency diagram corresponding to the maximum of the relative error in the energy along the integration interval.
In this case, the new scheme $\mathcal{A}_{18}$ is the most efficient when high accuracy is desired. There are initial conditions, however, for which 
RKN6$_{11}$ provides better results up to round-off.

Very similar results are obtained for the H\'enon-Heiles potential, and for this reason they are not shown here. From the previous experiments, we can conclude
that the new scheme  $\mathcal{A}_{19}$ outperforms the symplectic methods of order 4 and 6 from medium to high accuracy when the potential has a singularity,
whereas $\mathcal{A}_{17}$, $\mathcal{A}_{18}$ and $\mathcal{B}_{18}$ deliver the best results only at high accuracy for smooth potentials. To provide
further evidence to this class, we next consider a slightly more involved example.

\paragraph{Example 7: the restricted three body problem.}
In this case we have two bodies of masses $1 - \mu$ and $\mu$ in
circular rotation in a plane and a third body of negligible mass moving around in
the same plane. The equations of motion in a fixed coordinate system read \cite[p. 129]{hairer93sod}
\begin{equation}\label{eq.AREnstorf1}
\begin{aligned}
   & \ddot y_1 = y_1+2\dot y_2 - \mu' \, \frac{y_1+\mu}{D_1} - \mu \, \frac{y_1-\mu'}{D_2} \\
   &  \ddot y_2 = y_2-2\dot y_1 - \mu'\frac{y_2}{D_1} - \mu\frac{y_2}{D_2}, 
\end{aligned}   
\end{equation}
where $D_1=((y_1+\mu)^2+y_2^2)^{3/2}$, $D_2=((y_1-\mu')^2+y_2^2)^{3/2}$, and $\mu'= 1-\mu$. This system can be split as in \eqref{rkn.1b}-\eqref{rkn.2b}. 
Alternatively, in a rotating system the equations of motion become
\begin{equation}\label{eq.AREnstorf2}
\begin{aligned}
   & \ddot y_1 = \mu' \, \frac{a_1(t)-y_1}{D_1} + \mu \, \frac{b_1(t)-y_1}{D_2} \\
   & \ddot y_2 = \mu' \, \frac{a_2(t)-y_2}{D_1} + \mu \, \frac{b_2(t)-y_2}{D_2}, 
\end{aligned}      
\end{equation}
where now
\[
D_1=((y_1-a_1(t))^2+(y_2-a_2(t))^2)^{3/2}, \quad D_2=((y_1-b_1(t))^2+(y_2-b_2(t))^2)^{3/2}, 
\]
and the motion of the massive bodies is described by
\[
  a_1(t)=-\mu\cos(t), \quad
  a_2(t)=-\mu\sin(t); \quad
  b_1(t)=\mu'\cos(t), \quad
  b_2(t)=\mu'\sin(t). \quad
\]
We take, as in \cite{hairer93sod},  $\mu = 0.012277471$ and the following initial conditions in the rotating system:
\[
 y_1(0) = 0.994, \quad
 \dot y_1(0) = 0,  \quad
 y_2(0) = 0,  \quad
 \dot y_2(0) = 
-1.00758510637908252240.
\]
The resulting closed trajectory corresponds to the so-called Arenstorf orbit in the fixed coordinate system, with period
$T = 17.06521656015796255889$. 

\begin{figure}[!h]
  \begin{center}
      \includegraphics[width=10cm]{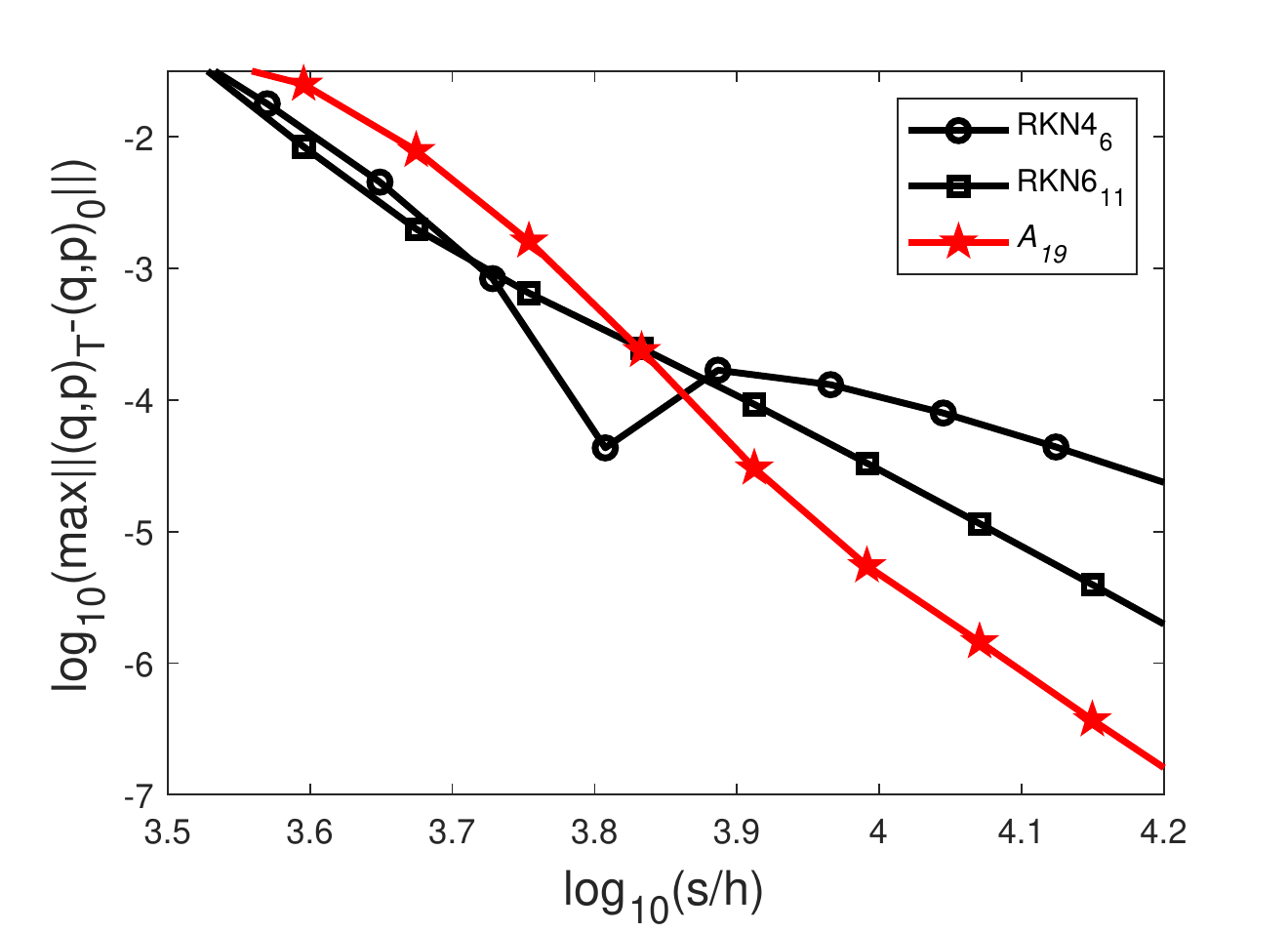}
    \caption{\small Error with respect to the initial conditions after one period, $T$, of the Arenstorf orbit versus the number of force evaluations for the 4th-, 6th- and 8th-order RKN splitting methods, RK4N$_6$ (circles), RKN6$_{11}$ (squares) and $\mathcal{A}_{19}$ (stars). }
    \label{fig:ARE_RKN}
  \end{center}
\end{figure}

In this case we integrate for one period with the RKN splitting methods of order 4 and 6, and the new 8th-order scheme $\mathcal{A}_{19}$. We
measure the error with respect to the initial conditions (taking into account that we are integrating in the rotating system) and display the
corresponding errors in Figure~\ref{fig:ARE_RKN}. Again, $\mathcal{A}_{19}$ is the most efficient scheme even for medium accuracies.

\section{Conclusions}

We have presented new RKN splitting methods of order 8 that show a better efficiency than the best existing symmetric compositions
of 2nd-order symmetric schemes on a variety of examples. We have thus answered in the affirmative the question formulated by \cite{mclachlan96aso} 
in 1996 and filled the existing gap in the classification of the most efficient splitting and composition methods \cite{blanes08sac,mclachlan02sm}. The
technical difficulties involved in the process have been overcome by applying standard techniques for solving nonlinear polynomial equations
 and free software on a personal computer. Whereas previous 8th-order RKN splitting methods require the evaluation of `modified potentials' or
 force-gradients \cite{omelyan02otc}, 
 the  schemes collected here only involve the evaluation of the force $g(y)$, just as compositions (\ref{eq.2.1.2}) and thus they should
 be considered as the natural option when one is interested in integrating the system (\ref{rkn.1}) with high precision and the evaluation of modified
 potentials is computationally expensive or not feasible.
 
 Both types of compositions (\ref{aba}) and (\ref{bab}) have been analyzed and different schemes with up to two free parameters have been constructed
 and tested on different numerical examples. These show that $\mathcal{A}_{18}$ and $\mathcal{B}_{18}$ provide better efficiencies when the force is
 derived from a smooth, singularity-free potential, whereas for problems involving singularities $\mathcal{A}_{19}$ exhibits the best results.
 As representatives of the first situation
 (i.e., singularity-free potentials), we have examined the simple pendulum, the H\'enon--Heiles potential and
 the quantum treatment of the P\"oschl--Teller potential. The second case, involving singularities, corresponds
 to the Kepler problem and the restricted planar three body problem. Moreover, the new schemes
 are more efficient than lower order RKN splitting methods for medium to high accuracies, and provide better results than extrapolation methods
 of order 8 even for relatively short time integrations.

\subsection*{Acknowledgements}
 This work has been supported by 
Ministerio de Ciencia e Innovaci\'on (Spain) through project PID2019-104927GB-C21, MCIN/AEI/10.13039/501100011033.
A.E.-T. has been additionally funded by the predoctoral contract BES-2017-079697 (Spain).
A.E.-T. would like to thank Ander Murua and Joseba Makazaga for their help in implementing their continuation algorithms and the UPV-EHU
for its hospitality.


\begin{thebibliography}{10}

\bibitem{makazaga00ac}
\texttt{https://github.com/jmakazaga/arc-continuation}


\bibitem{alberdi19aab}
{\sc E.~Alberdi, M.~Anto{\~n}ana, J.~Makazaga, and A.~Murua}, {\em An algorithm
  based on continuation techniques for minimization problems with highly
  non-linear equality constraints}, Tech. Rep. arXiv:1909.07263, 2019.

\bibitem{arnold89mmo}
{\sc V.~Arnold}, {\em Mathematical {M}ethods of {C}lassical {M}echanics},
  Springer-Verlag, {S}econd~ed., 1989.

\bibitem{blanes16aci}
{\sc S.~Blanes and F.~Casas}, {\em A {C}oncise {I}ntroduction to {G}eometric
  {N}umerical {I}ntegration}, {CRC} Press, 2016.

\bibitem{blanes08sac}
{\sc S.~Blanes, F.~Casas, and A.~Murua}, {\em Splitting and composition methods
  in the numerical integration of differential equations}, Bol. Soc. Esp. Mat.
  Apl., 45 (2008), pp.~89--145.

\bibitem{blanes01hor}
{\sc S.~Blanes, F.~Casas, and J.~Ros}, {\em High-order
  {R}unge--{K}utta--{N}ystr\"om geometric methods with processing}, Appl.
  Numer. Math., 39 (2001), pp.~245--259.

\bibitem{blanes01nfo}
{\sc S.~Blanes, F.~Casas, and J.~Ros}, {\em New families of
  symplectic {R}unge--{K}utta--{N}ystr\"om integration methods}, in Numerical
  Analysis and its Applications, LNCS 1988, Springer, 2001, pp.~102--109.

\bibitem{blanes02psp}
{\sc S.~Blanes and P.~Moan}, {\em Practical symplectic partitioned
  {R}unge--{K}utta and {R}unge--{K}utta--{N}ystr\"om methods}, J. Comput. Appl.
  Math., 142 (2002), pp.~313--330.

\bibitem{calvo93hos}
{\sc M.~Calvo and J.~Sanz-Serna}, {\em High-order symplectic
  {R}unge--{K}utta--{N}ystr\"om methods}, SIAM J. Sci. Comput., 14 (1993),
  pp.~1237--1252.

\bibitem{flugge71pqm}
{\sc S.~Fl\"ugge}, {\em Practical Quantum Mechanics}, Springer, 1971.

\bibitem{hairer06gni}
{\sc E.~Hairer, C.~Lubich, and G.~Wanner}, {\em Geometric {N}umerical
  {I}ntegration. {S}tructure-{P}reserving {A}lgorithms for {O}rdinary
  {D}ifferential {E}quations}, Springer-Verlag, {S}econd~ed., 2006.

\bibitem{hairer93sod}
{\sc E.~Hairer, S.~N{\o}rsett, and G.~Wanner}, {\em Solving {O}rdinary
  {D}ifferential {E}quations {I}, {N}onstiff {P}roblems}, Springer-Verlag,
  {S}econd revised~ed., 1993.

\bibitem{numpy}
{\sc C.~R. Harris, K.~J. Millman, S.~J. van~der Walt, R.~Gommers, P.~Virtanen,
  D.~Cournapeau, E.~Wieser, J.~Taylor, S.~Berg, N.~J. Smith, R.~Kern, M.~Picus,
  S.~Hoyer, M.~H. van Kerkwijk, M.~Brett, A.~Haldane, J.~Fern\'andez~del R\'{\i}o,
  M.~Wiebe, P.~Peterson, P.~G\'erard-Marchant, K.~Sheppard, T.~Reddy,
  W.~Weckesser, H.~Abbasi, C.~Gohlke, and T.~E. Oliphant}, {\em Array
  programming with NumPy}, Nature, 585 (2020), p.~357--362.

\bibitem{henon64tao}
{\sc M.~H\'enon and C.~Heiles}, {\em The applicability of the third integral of
  motion: some numerical experiments}, Astron. J., 69 (1964), pp.~73--79.

\bibitem{kahan97ccf}
{\sc W.~Kahan and R.~Li}, {\em Composition constants for raising the order of
  unconventional schemes for ordinary differential equations}, Math. Comput.,
  66 (1997), pp.~1089--1099.

\bibitem{krafft2005debian}
{\sc M.~F. Krafft}, {\em The Debian System: Concepts and Techniques}, No Starch
  Press, 2005.

\bibitem{lubich08fqt}
{\sc C.~Lubich}, {\em From {Q}uantum to {C}lassical {M}olecular {D}ynamics:
  {R}educed {M}odels and {N}umerical {A}nalysis}, European Mathematical
  Society, 2008.


\bibitem{mclachlan19tla}
{\sc R.~McLachlan and A.~Murua}, {\em The {L}ie algebra of classical
  mechanics}, J. Comput. Dyn., 6 (2019), pp.~198--213.

\bibitem{mclachlan02sm}
{\sc R.~McLachlan and R.~Quispel}, {\em Splitting methods}, Acta Numerica, 11
  (2002), pp.~341--434.

\bibitem{mclachlan03tae}
{\sc R.~McLachlan and B.~Ryland}, {\em The algebraic entropy of classical
  mechanics}, J. Math. Phys., 44 (2003), pp.~3071--3087.

\bibitem{mclachlan96aso}
{\sc R.~McLachlan and C.~Scovel}, {\em A survey of open problems in symplectic
  integration}, in Integration {A}lgorithms and {C}lassical {M}echanics,
  J.~Marsden, G.~Patrick, and W.~Shadwick, eds., Fields Institute
  Communications, American Mathematical Society, 1996, pp.~151--180.

\bibitem{minpack}
{\sc J.~J. Mor{\'e}, B.~S. Garbow, and K.~E. Hillstrom}, {\em User guide for
  {MINPACK}-1}, tech. rep., CM-P00068642, 1980.

\bibitem{okunbor94oeo}
{\sc D.~Okunbor and E.~Lu}, {\em Eighth-order explicit symplectic
  {R}unge--{K}utta--{N}ystr\"om integrators}, Tech. Rep. CSC 94-21, University
  of Missouri-Rolla, 1994.

\bibitem{omelyan02otc}
{\sc I.~Omelyan, I.~Mryglod, and R.~Folk}, {\em On the construction of high
  order force gradient algorithms for integration of motion in classical and
  quantum systems}, Phys. Rev. E, 66 (2002), p.~026701.

\bibitem{sofroniou05dos}
{\sc M.~Sofroniou and G.~Spaletta}, {\em Derivation of symmetric composition
  constants for symmetric integrators}, Optim. Method. Softw., 20 (2005),
  pp.~597--613.

\bibitem{trefethen00smi}
{\sc L.~Trefethen}, {\em Spectral {M}ethods in {MATLAB}}, {SIAM}, 2000.

\bibitem{python3rm}
{\sc G.~Van~Rossum and F.~L. Drake}, {\em Python 3 Reference Manual},
  CreateSpace, Scotts Valley, CA, 2009.

\bibitem{varadarajan84lgl}
{\sc V.~Varadarajan}, {\em Lie {G}roups, {L}ie {A}lgebras, and {T}heir
  {R}epresentations}, Springer-Verlag, 1984.

\bibitem{scipy}
{\sc P.~Virtanen, R.~Gommers, T.~E. Oliphant, M.~Haberland, T.~Reddy,
  D.~Cournapeau, E.~Burovski, P.~Peterson, W.~Weckesser, J.~Bright, S.~J. {van
  der Walt}, M.~Brett, J.~Wilson, K.~J. Millman, N.~Mayorov, A.~R.~J. Nelson,
  E.~Jones, R.~Kern, E.~Larson, C.~J. Carey, {\.I}.~Polat, Y.~Feng, E.~W.
  Moore, J.~{VanderPlas}, D.~Laxalde, J.~Perktold, R.~Cimrman, I.~Henriksen,
  E.~A. Quintero, C.~R. Harris, A.~M. Archibald, A.~H. Ribeiro, F.~Pedregosa,
  P.~{van Mulbregt}, and {SciPy 1.0 Contributors}}, {\em {{SciPy} 1.0:
  Fundamental algorithms for scientific computing in Python}}, Nature Methods,
  17 (2020), pp.~261--272.

\end{thebibliography}

\end{document}